\documentclass[12pt, leqno, amstex, latexsym, amscd, amssymb]{article}
\usepackage{amsmath,amssymb,amscd,latexsym}

\def\noproof{{\unskip\nobreak\hfill\penalty50\hskip2em\hbox{}%
     \nobreak\hfill$\Box$\parfillskip=0pt%
     \finalhyphendemerits=0\par}}
\def\enddemo{\ifmmode\eqno\Box\else\noproof\vskip0.8truecm\fi}
\newcommand{\noi}{\noindent}

\newtheorem{theorem}{Theorem}[section]
\newtheorem{definition}[theorem]{Definition}
\newtheorem{conjecture}[theorem]{Conjecture}

\newtheorem{proposition}[theorem]{Proposition}
\newtheorem{exam}[theorem]{Examples}

\newtheorem{coro}[theorem]{Corollary}

\newtheorem{lemma}[theorem]{Lemma}
\newtheorem{remark}[theorem]{Remark}
\newtheorem{prop}[theorem]{Proposition}
\newtheorem{assumption}[theorem]{Assumption}

\newcommand{\lra}{\longrightarrow}
\newcommand{\Gal}{\mbox{Gal}}
\newcommand{\sym}{\mbox{Sym}}
\newcommand{\red}{\mbox{red}}

\newcommand{\mat}[4]{\left( \begin{array}{cc} {#1} & {#2} \\ {#3} & {#4}\end{array} \right)}
\newcommand{\edges}{{\stackrel{\rightarrow}{{\cal E}}}}
\newcommand{\vertices}{{\cal V}}

\newcommand{\Image}{\mbox{Im}}

\newcommand{\wh}{\widehat}
\newcommand{\wt}{\widetilde}
\newcommand{\pa}{\partial}

\newcommand{\ord}{\mbox{ord}}
\newcommand{\Pic}{\mbox{Pic}}
\newcommand{\NS}{\mbox{NS}}

\newcommand{\spec}{\mbox{Spec}}
\newcommand{\spf}{\mbox{Spf}}
\newcommand{\End}{\mbox{End}}
\newcommand{\Aut}{\mbox{Aut}}
\newcommand{\ad}{\mbox{ad}}
\newcommand{\Norm}{\mbox{Nrd}}
\newcommand{\trace}{\mbox{trace}}

\newcommand{\PGL}{{\rm \bold{PGL}}}
\newcommand{\GL}{{\rm \bold{GL}}}

\newcommand{\SL}{{\rm {\bold{SL}}}}
\newcommand{\GG}{{\rm {\bold{G}}}}

\newcommand{\bu}{{\bullet}}

\newcommand{\Ext}{\mbox{Ext}}
\newcommand{\Ind}{\mbox{Ind}}
\newcommand{\Maps}{\mbox{Maps}}
\newcommand{\Cone}{\mbox{Cone}}
\newcommand{\Res}{\mbox{Res}}
\newcommand{\Ker}{\mbox{Ker}}
\newcommand{\Hom}{\mbox{Hom}}
\newcommand{\Corr}{\mbox{Corr}}
\newcommand{\emb}{\mbox{emb}}
\newcommand{\Rep}{{\rm \mbox{\underline{Rep}}}}
\newcommand{\seq}{\subseteq}
\newcommand{\cont}{{\scriptsize{\mbox{cont}}}}
\newcommand{\sRep}[1]{{\scriptsize{\mbox{\underline{Rep}$_{ss}(G_{#1})$}}}}

\newcommand{\et}{{\scriptsize{\mbox{\'et}}}}
\newcommand{\st}{{\scriptsize{\mbox{st}}}}
\newcommand{\cris}{{\scriptsize{\mbox{cris}}}}
\newcommand{\DR}{{\scriptsize{\mbox{DR}}}}
\newcommand{\har}{{\scriptsize{\mbox{har}}}}

\newcommand{\Abar}{{\overline{A}}}

\newcommand{\Kbar}{{\overline{K}}}
\newcommand{\Ubar}{{\overline{U}}}
\newcommand{\Xbar}{{\overline{X}}}
\newcommand{\Zbar}{{\overline{Z}}}

\newcommand{\Gbar}{{\overline{G}}}
\newcommand{\Fpbar}{{\overline{\F}_p}}
\newcommand{\kbar}{{\overline{k}}}

\newcommand{\Qbar}{{\overline{\Q}}}

\newcommand{\Pbar}{{\overline{P}}}
\newcommand{\Psibar}{{\overline{\Psi}}}

\newcommand{\cK}{{\cal K}}
\newcommand{\cL}{{\cal L}}
\newcommand{\cM}{{\cal M}}
\newcommand{\cH}{{\cal H}}

\newcommand{\cP}{{\cal P}}
\newcommand{\cT}{{\cal T}}
\newcommand{\PP}{{\Bbb P}}

\newcommand{\cO}{{\cal O}}
\newcommand{\cA}{{\cal A}}
\newcommand{\cE}{{\cal E}}
\newcommand{\cV}{{\cal V}}

\newcommand{\cF}{{\cal F}}

\newcommand{\cG}{{\cal G}}

\newcommand{\cB}{{\cal B}}
\newcommand{\cR}{{\cal R}}
\newcommand{\cU}{{\cal U}}
\newcommand{\cW}{{\cal W}}
\newcommand{\cX}{{\cal X}}

\newcommand{\TT}{{\Bbb T}}
\newcommand{\Z}{{\Bbb Z}}
\newcommand{\Q}{{\Bbb Q}}
\newcommand{\Ql}{{\Q_{\ell}}}
\newcommand{\C}{{\Bbb C}}

\newcommand{\F}{{\Bbb F}}
\newcommand{\BGG}{{\Bbb G}}

\newcommand{\bL}{{\Bbb L}}
\newcommand{\bT}{{\bf T}}

\newcommand{\bH}{{\bf H}}
\newcommand{\DD}{{\Bbb D}}

\newcommand{\Qur}{{\Q^{ur}_p}}
\newcommand{\Zur}{{\Z^{ur}_p}}

\newcommand{\fz}{{\frak{z}}}

\newcommand{\fa}{{\frak{a}}}

\newcommand{\fX}{{\frak{X}}}
\newcommand{\fY}{{\frak{Y}}}

\newcommand{\s}{{\sigma}}

\newcommand{\p}{{\phi}}
\newcommand{\G}{{\Gamma}}
\newcommand{\la}{{\lambda}}
\newcommand{\La}{{\Lambda}}
\newcommand{\De}{{\Delta}}
\newcommand{\ep}{{\epsilon}}
\newcommand{\io}{{\iota}}
\newcommand{\ga}{{\gamma}}
\newcommand{\ka}{{\kappa}}

\newcommand{\Hx}{{H^1_{\DR}(X_{\G}, \cV)}}
\newcommand{\Hxn}{{H^1_{\DR}(X_{\G}, \cV_n)}}
\newcommand{\Hu}{{H^1_{\DR}(U, \cV)}}
\newcommand{\Hcu}{{H^1_{\DR,c}(U, \cV)}}
\newcommand{\Hun}{{H^1_{\DR}(U, \cV_n)}}
\newcommand{\Hcun}{{H^1_{\DR,c}(U, \cV_n)}}

\begin{document}

\title{Derivatives of $p$-adic $L$-functions, Heegner cycles and monodromy
modules attached to modular forms}
\author{By Adrian Iovita and Michael Spiess}
\date{}
\maketitle
\tableofcontents

\section{Introduction}
\label{sec:introduction}

Let $f(z)$ be an elliptic newform of even weight $k=n+2\geq 4$ and level $N$ 
and let $p$ be a prime number. Let us assume that $f(z)$ corresponds to a 
modular form  on a Shimura curve via the theory of Jaquet-Langlands
and let us fix a quadratic imaginary field $K$. Throughout the introduction
we will assume for simplicity that $N$ is square free with an even number of 
prime factors, $p|N$, all prime divisors of $N$ are inert in $K$ and
finally that $K$ has class number 1 (see section \ref{sec:dieumod} and \ref{sec:heegner} for a more general set-up).  

In such a situation in \cite{bdis} we have defined the anticyclotomic $p$-adic $L$-function of $f(z)$ over $K$, $L_p(f/K, s)$ and have proved that it vanishes at the central critical point $s=k/2$. The main purpose of this paper is to prove that the derivative of $L_p(f/K, s)$ at $s=k/2$ can be interpreted as the 
$e_f$-component ($e_f$ is the idempotent in the Hecke algebra corresponding to 
$f(z)$) of the image of a Heegner cycle under a $p$-adic Abel-Jacobi map. The complex $L$-function of $f(z)$ over $K$ also vanishes at $s=k/2$ and its derivative can be expressed in terms of the height of a Heegner cycle. Therefore we consider our result as a $p$-adic Gross-Zagier type formula though the reader should be aware that no $p$-adic heights are involved (a short explanation of how our result fits into the broader picture of the Bloch-Beilinson conjectures is given in Remark \ref{remark:selmer}). 

With some slight modifications our argument would also give a similar result for the weight $k=2$. If $f(z)$  corresponds to an elliptic curve $E/\Q$ this is then a reformulation of the main result of \cite{bertolini_darmon2}. We also remark that a $p$-adic Gross-Zagier type formula involving $p$-adic heights has been obtained by Nekovar \cite{nekovar3}. However the situation he deals with is different; he considers a two-variable $p$-adic $L$-function and the case where $p$ does not devide the level and $f(z)$ is {\it ordinary} at $p$. 

In establishing our formula we make extensive use of $p$-adic Hodge theory.
Namely, if $G_\Q\to \GL(V_p(f))$ is the Galois representation attached to
$f(z)$ as in \cite{deligne1}, its restriction to a decomposition group at $p$
 is a semistable representation. Our proof relies on the fact that we are able to explicitly describe the semistable Dieudonn\'e module attached to $V_p(f)$ (via Fontaine's theory) in terms of $p$-adic integration. We think that this description is interesting in itself and as a separate application we use it to prove that the $\cL_p$-invariants of Fontaine-Mazur and Teitelbaum attached to $f(z)$ are equal.\\
 
Now we describe our results in more detail. Let $X/\Q$ be the Shimura curve 
which classifies abelian surfaces with an action of a maximal order in the 
quaternion algebra over $\Q$ of discriminant $N$. The Jacquet-Langlands 
theory associates to $f(z)$ a Hecke eigenform on $X$ which by abuse of 
notation we also denote by $f(z)$. Let $F_f$ be the finite extension of 
$\Q$ generated by the Hecke eigenvalues of $f(z)$. The Galois representation 
$V_p(f)$ attached to $f(z)$ (two-dimensional as a $F_f\otimes \Q_p$-module) 
occurs in the $p$-adic cohomology of a certain Kuga-Sato variety. It can be 
also realized as a direct summand of the $(n+1)$-th $p$-adic cohomology of 
the $m \colon = \frac{n}{2}$-self product $\cA^m$ of the universal abelian 
surface $\cA$ over $X$ (in the introduction we neglect the fact that $X$ is 
only a coarse moduli space; for precise statements we refer to section 
\ref{sec:dieumod}). In fact the Galois representation attached to the 
whole space of modular forms of weight $k$ which we denote by $H_p(\cM_n)$ 
(since it is the $p$-adic realisation of a certain motive $\cM_n$) can be 
identified with the first \'etale cohomology of $\Xbar = X \otimes_{\Q} 
\Qbar$ with coefficients in a $p$-adic local system $\bL_n$ which occurs in 
a certain tensor power of the first relative $p$-adic cohomology of $\cA$ 
over $X$. As a representation of $G_{\Q_p}$, $H_p(\cM_n)$ is semisimple. 
We can apply Fontaine's theory which associates to local $p$-adic Galois 
representations simpler objects -- filtered $(\p,N)$-modules -- which still 
encode all the information. When restricting the Galois action to the inertia 
subgroup $I_p\seq G_{\Q_p}$ and applying Fontaine's functor $D_{\st}$ we 
obtain a filtered $(\p, N)$-module $D_{\st}(H_p(\cM_n))$ the 
{\it semistable Dieudonn{\'e} module} of $H_p(\cM_n)$. That is 
$D_{\st}(H_p(\cM_n))$ is a finite dimensional vector space over the 
completion $\Qur$ of the maximal unramified extension of $\Q_p$ endowed with 
Frobenius $\p$ and monodromy operator $N$. In section \ref{sec:dieumod} we give a concrete description of $D_{\st}(H_p(\cM_n))$ by using the 
\v Cerednik-Drinfeld uniformisation of $X_{\Qur} = X\otimes_{\Q} \Qur$ by 
the $p$-adic upper half plane $\cH_p$. More precisely we can identify 
$X_{\Qur}$ with the Mumford curve $X_{\G} = \G \backslash \cH_p$ (the group 
$\G\seq \GL_2(\Q_p)$ is given in terms of the definite quaternion algebra 
over $\Q$ of discriminant $N/p$). We show that $D_{\st}(H_p(\cM_n))$ is equal 
to $H^1_{\DR}(X_{\G}, \cV_n)$ where $\cV_n$ is a filtered $F$-isocrystal on 
$X_{\G}$ associated to the $\GL_2(\Q_p)$-representation $V_n \colon = n-\mbox{th symmetric power of $\Q_p^2$}$. The Frobenius operator is given in terms of Coleman integration of $1$-forms on $\cH_p$ and the monodromy operator in terms of Schneider's residue map $\omega \mapsto \Res_e(\omega)$ which associates to a $1$-form a harmonic cocyle on the Bruhat-Tits tree (see Theorem \ref{theorem:dieumodform}). We also have a result for open subschemes of $X$ which is needed for the computation of the $p$-adic Abel-Jacobi image of Heegner cycles.

The main ingredients in the proof of Theorem \ref{theorem:dieumodform} 
are a comparison theorem between $p$-adic \'etale cohomology of semistable 
curves and log-crystalline cohomology (with coefficients) in 
 \cite{faltings2} and the description of the log-crystalline cohomology groups 
in \cite{coleman_iovita}.  These facts are reviewed in section 
\ref{sec:isocrystals}. In section \ref{sec:gltwo} we attach to a 
$\Q_p$-rational representation of $\GL_2$ a  filtered $F$-isocrystal on 
$X_{\G}$.  A key result is then that the local system corresponding to the 
filtered $F$-isocrystal attached to $V_n$ is the sheaf $\bL_n$ (see Lemma \ref{lemma:h1dr}).\\

The application to $\cL$-invariants is given in section \ref{sec:linv}. Recall that the $\cL$-invariant for a weight two modular form $f(z)$ corresponding to an elliptic curve $E/\Q$ with split multiplicative reduction at $p$ is defined as the ratio $\cL(f) = \log_p(q_E)/\ord_p(q_E)$ where $q_E$ is the Tate period of $E$ over $\Q_p$. In the higher weight case three possible definitions of $\cL(f)$ have been given (by Teitelbaum \cite{teitelbaum}, Coleman \cite{cole} and Fontaine-Mazur \cite{mazur}). The first and second are defined in terms of 
Coleman integration and residues on the Shimura curve and modular curve respectively whereas the last in terms of the semistable Dieudonn{\'e} module $D_{\st}(V_p(f))$. As an application of Theorem \ref{theorem:dieumodform} we will show that the first and last are equal. We remark that the work of \cite{coleman_iovita} also establishes the equality of the Coleman and Fontaine-Mazur $\cL$-invariant (as explained in \cite{cole}) so that all three $\cL$-invariants are the same. 

In section \ref{sec:lfunc} we apply our  results to  Heegner cycles and 
$p$-adic $L$-functions. In \cite{bdis} we proved the following formula 
describing the derivative of \linebreak
$L_p(f/K, s)$ at $s=k/2$ as a $p$-adic integral
\begin{equation}
\label{eqn:intromain}
L_p'(f/K,k/2) = \int^{z_0}_{\bar{z}_0}\,\, f_p(z) P(z)^m dz.
\end{equation}
Here $f_p(z)$ is the rigid analytic modular form on $\cH_p$ for $\G$ 
corresponding to $f(z)$, $P(z)$ a certain polynomial of degree $2$ and 
$z_0$ is a point on $\cH_p$ lying over a Heegner point on $X$.

In section \ref{sec:lfunc} we show  that the right hand side of (\ref{eqn:intromain}) can be interpreted as the ($e_f$-component) of the image of a Heegner cycle on $\cA^m$ (defined in section \ref{sec:heegner}) under a $p$-adic Abel-Jacobi map (the latter will be defined in section \ref{sec:pabeljac}). We briefly explain the main steps of the proof. Firstly a Heegner cycle defines an extension of $D_{\st}(H_p(\cM_n))[-(m+1)]$ by the trivial $(\p,N)$-module $\Qur$. The group of extension classes $\Ext^1(\Qur, D_{\st}(H_p(\cM_n)[-(m+1)])$ can be identified with the dual of the space $M_k(\G)$ of weight $k$ modular forms on $\cH_p$ for $\G$. It is then shown that under this isomorphism the extension class corresponding to the Heegner cycle is the functional 
\[
f_p(z) \mapsto \int^{z_0}_{\bar{z}_0}\,\, f_p(z) P(z)^m dz.
\]
This last step uses a formula expressing the cup-product of 1-forms on $X_{\G}$ in terms of the periods (this formula is a generalisation of a result of E. de Shalit in \cite{shalit2}).\\

\bigskip\noindent
{\bf Aknowledgements}  We would like to thank Massimo Bertolini and 
Henri Darmon for sharing with us their ideas on anticyclotomic $p$-adic 
$L$-functions. We are also grateful to Amnon Besser and Uwe Jannsen for 
helpful discussions. Both authors have visited the SFB-478  M\"unster 
during the writing of this paper and thank this institution for its 
support and hospitatlity. This work was supported in part by EPSRC Grant 
GR/M95615. During the work on this project the first author was also 
partially supported by an NSF grant.

\bigskip
\noi {\bf Notation.}

For a field $K$ we let $\Kbar$ be an algebraic closure and $G_K = \Gal(\Kbar/K)$ the absolute Galois group.

Let $L/K$ be a field extension. For a $K$-module $M$ we set $M_L = M\otimes_K L$. Also for a $K$-scheme $X$ we write $X_L = X \otimes_K L$.

For a $\Q$-module $M$ we write $M_p = M \otimes_{\Q} \Q_p$.

We denote by $\C_p$ the completion of the algebraic closure of $\Q_p$. We let $\Qur\seq \C_p$ be the closure of the maximal unramified extension of $\Q_p$ and $\Zur$ its valuation ring. For a positive integer $\mu$ we let $\Q_{p^{\mu}}\seq \Qur$ be the unramified extension of degree $\mu$ of $\Q_p$ and $\Z_{p^{\mu}}$ its ring of integers.

\section{Filtered Frobenius monodromy modules}
\label{sec:phin}

Let $K$ be a field of characteristic 0 which is complete with respect to a discrete valuation and has a perfect residue field $\ka$ of characteristic $p> 0$. Let $K_0\seq K$ be the maximal subfield of $K$ of absolute ramification index $1$, i.e.\ $K_0$ is the quotient field of the ring of Witt vectors of $\ka$.

Let $\s: K_0 \to K_0$ be the absolute Frobenius automorphism. In this section we deal with $\s$-linear algebra. Specifically we recall (form  \cite{font2}) the relation between semistable $p$-adic representations of $G_K$ and filtered Frobenius monodromy modules (filtered $(\p, N)$-modules for short) and then give a description of certain $\Ext$-groups in the category of (filtered $(\p, N)$-modules) which will be used in our discussion of the $p$-adic Abel-Jacobi map in section \ref{sec:pabeljac}. We also introduce the notion of monodromy modules and their $\cL$-invariants which is needed for the definition of the Fontaine-Mazur $\cL$-invariant of modular forms.

A filtered $(\p, N)$-module consists of the following data:\\
1. A finite dimensional $K_0$-vector space $D$ with an exhaustive and separated decreasing filtration $F^i D_K$  on $D_K$ called Hodge filtration.\\
2. A $\s$-linear automorphism $\p = \p_D: D\lra D$ (the Frobenius of $D$).\\
3. A $K$-linear endomorphism $N = N_D: D \lra D$ (the monodromy operator) such that $N \p = p \p N$.

A filtered Frobenius module is a $(\p, N)$-module with trivial monodromy operator. Morphisms of filtered $(\p, N)$-modules are $K_0$-linear maps which respect the Frobenii, the filtrations and the monodromy operators. The category of filtered $(\p, N)$-modules $MF_{K}(\p, N)$ is an additive tensor category admitting kernels and cokernels. There is also the notion of a short exact sequence of filtered $(\p, N)$-modules. We consider $K$ as a filtered $(\p, N)$-module by setting $\p_{K_0} = \s$, $N = 0$ and  $F^i K = K$ (resp.\ $=0$) for $i\leq 0$ (resp.\ $i>0$). For $D\in Ob(MF_{K}(\p, N))$ its $i$-fold twist $D[i]$ is defined as $D[i] = D$ (as vector spaces), $\p_{D[i]} = p^i \p_D$, $F^j D[i]_K = F^{j-i} D_K$ and $N_{D[i]} = N_D$.

The correspondence between filtered $(\p, N)$-modules and semistable Galois representations is given in terms of Fontaine's ring $B_{\st}$ (defined in \cite{font1}). It is a topological $K_0$-algebra whose construction depends on chosing a branch of the $p$-adic logarithm. It is equipped with the following structure:\\
1. A continous action of $G_{K}$ such that $B_{st}^{G_{K}} = K_0$.\\
2. A $G_{K}$-equivariant embedding $K_0^{ur} \lra B_{st}$ of the maximal unramified extension $K_0^{ur}$ of $K_0$.\\
3. A $\s$-linear continous automorphism $\p: B_{st}\lra B_{st}$ commuting with the $G_{K}$-action.\\
4. An exhaustive and separated decreasing filtration $F^i$ on $(B_{st})_K$ which is stable under the $G_{K}$-action.\\
5. A $K_0$-linear operator $N: B_{st} \lra B_{st}$ such that $N \p = p \p N$. 

Let $\Rep(G_{K})$ be the category of $p$-adic representations of $G_{K}$, i.e.\ finite dimensional $\Q_p$-vector spaces with a continuous linear $G_{K}$-action. It is an abelian tensor category with twists given by tensoring with appropriate powers of the Tate representation $\Q_p(1) = (\underset{\longrightarrow}{\lim} \mu_{p^n})\otimes_{\Z_p} \Q_p$. For a $p$-adic representation $V$ of $G_{K}$ one defines $D_{\st}(V)\colon =(V\otimes B_{st})^{G_{K}}$. It inherits from $B_{st}$ the structure of a filtered $(\p, N)$-module. We have $D_{st}(V(n)) = D_{st}(V)[-n]$. A $p$-adic representation $V$ of $G_{K}$ is called semistable if the canonical (injective) map 
\[
\alpha: D_{\st}(V)\otimes_{K_0} B_{st} \lra V\otimes_{\Q_p} B_{st} 
\]
is bijective. The full subcategory $\Rep_{ss}(G_{K})$ of semistable representations is an abelian tensor category. 

Conversely given a filtered $(\p, N)$-module $D$ the module $D\otimes_{K_0} B_{st}$ has a natural structure as a filtered $(\p, N)$-module and one defines $V_{st}(D)\colon = \Hom_{MF_{K}(\p, N)}(K, D \otimes_{K_0} B_{st})$. It is a $p$-adic representation of $G_{K}$. The module $D$ is called admissible if it is isomorphic to $D_{st}(V)$ for some semistable representation $V$ of $G_{K}$. The full subcategory $MF^{ad}_{K}(\p, N)$ of $MF_{K}(\p, N)$ of admissible filtered $(\p, N)$-modules is an abelian tensor category such that exact sequences remain exact in $MF_{K}(\p, N)$. Moreover the restriction of the functors $D_{st}$ and $V_{st}$ to $\Rep_{ss}(G_{K})$ and $MF^{ad}_{K}(\p, N)$ are mutually quasi-inverse $\otimes$-equivalences.\\

Any finite dimensional $K_0$-vector space $D$ equipped with a bijective $\s$-linear endomorphism $\p: D \to D$ admits a canonical {\it slope decomposition} (see \cite{zink})
\begin{equation}
\label{eqn:slope}
D  = \bigoplus_{\lambda\in \Q} D_{\lambda}.
\end{equation}
For a rational number $\la = \frac{r}{s}, r,s\in \Z, s>0$ the subspace $D_{\lambda}$ of $D$ (called the isotypical component of $D$ of slope $\lambda$) is the largest subspace of $D$ which has an $\cO_{K_0}$-stable lattice $M$ with $\p^s(M) = p^r M$. The {\it slopes} of $D$ are the rational numbers $\la$ for which $D_{\la} \neq 0$. The pair $(D,\p)$ is called isotypical of slope $\la_0$ if $\la_0$ is the only slope. If $D$ is a filtered $(\p, N)$-module then $N(D_{\lambda}) \seq D_{\lambda -1}$ for all $\lambda\in \Q$.

\begin{lemma} 
\label{lemma:bkexp}
Let $D$ be a filtered $(\p, N)$-module, $n$ be an integer and assume that $N$ induces an isomorphism between the isotypical components $D_n$ and $D_{n-1}$ . Then there is a canonical isomorphism 
\begin{equation}
\label{eqn:iso}
{\rm \Ext}^1_{MF_K(\p, N)}(K[n], D) \cong D /F^n.
\end{equation}
\end{lemma}

{\em Proof.} Firstly we remark that though $MF_K(\p, N)$ is not an abelian category a group structure on $\Ext^1_{MF_K(\p, N)}(K[n], D)$ can be defined in the usual way. For an extension 
\begin{equation}
\label{eqn:ext}
\begin{CD}
e:\quad 0 @>>> D @> f >> E @> r >> K[n] @>>> 0 
\end{CD}
\end{equation}
we obtain a diagram with exact rows
\[
\begin{CD}
0 @>>> D_n @>>> E_n @>>> K_0 @>>> 0\\
@.@VVV@VVV@VVV\\
0 @>>> D_{n-1} @>>> E_{n-1} @>>> 0 @>>> 0
\end{CD}
\]
where the vertical arrows are induced by the respective monodromy operators. By assumption the first vertical map is an isomorphism. Hence the upper row admits a canonical splitting that is the map $r: E \to K[n]$ has a canonical splitting $s: K[n] \to E$ (compatible with Frobenius and monodromy operators but not necessarily with filtrations).

We define the map (\ref{eqn:iso}) by sending the class of (\ref{eqn:ext}) to $s(1) +F^n \in E_K/ F^n \cong D_K /F^n$. A simple computation shows that (\ref{eqn:iso}) is an isomorphism.
\enddemo

We finish this section with a discussion of monodromy modules and its $\cL$-invariants. These notions were introduced in \cite{mazur} (in fact we work here with a slightly different definition; see Remark \ref{remark:Lbasechange} below). Let $\bT$ be a finite dimensional commutative semisimple $\Q_p$-algebra. For simplicity we assume from now on that $K = K_0$.

\begin{definition}
\label{definition:monodmod} Let $D$ be a $\bT$-object in $MF_K(\p, N)$, i.e.\ an object together with a $\Q_p$-algebra homomorphism of $\bT \to \End_{MF_K(\p, N)}(D)$. $D$ is called a (two-dimensional) monodromy $\bT$-module if the following conditions hold:\\
(i) $D$ is a free $\bT_K$-module of rank $2$.\\
(ii) The sequence $D \stackrel{N}{\lra} D \stackrel{N}{\lra} D$ is exact.\\
(iii) There exists an integer $j_0$ such that $F^{j_0} D$ is a free $\bT_K$-submodule of rank $1$ and $F^{j_0} D \cap \Ker(N_D) = 0$.
\end{definition}

Note that condition (i) and (ii) imply that $\Ker(N_D)$ is a free $\bT_K$-submodule of rank $1$. Hence (iii) implies that $F^{j_0} D \oplus \Ker(N_D) = D$.

\begin{lemma}
\label{lemma:phidecomp}
Let $D$ be a monodromy $\bT$-module. Then there exists a decomposition $D= D^{(1)}\oplus D^{(2)}$ where $D^{(1)}, D^{(2)}$ are $\p$-stable free rank-1 $\bT_K$-submodules such that $N: D \to D$ induces an isomorphism $N\mid_{D^{(2)}}:D^{(2)} \to D^{(1)}$. Moreover the decomposition is uniquely determined by theses properties.
\end{lemma}

{\em Proof}: We may assume that $\bT$ is a field. We claim that $D^{(1)} \colon = \Ker(N_D)$ is isotypical. For this it is enough to consider the case where $\ka$ is algebraically closed. Let $\la\in \Q$ with $D^{(1)}_{\la} \neq 0$. By a theorem of Dieudonn{\'e}-Manin the algebra $\End(D^{(1)}_{\la}, \phi)$ of the Frobenius isocrystal $(D^{(1)}_{\la}, \phi)$ is a central simple $\Q_p$-algebra of dimension $(\dim_K D^{(1)}_{\la})^2$. Since $\bT\seq \End(D^{(1)}_{\la})$ we deduce that $\dim_K D^{(1)}_{\la} \geq \dim_{\Q_p} \bT = \dim_K D^{(1)}$, hence $D^{(1)}_{\la} = D^{(1)}$ as claimed.

Property (ii) implies now that $D = D_{\la} \oplus D_{\la+1}$, $D^{(1)} = D_{\la}$ and $N: D_{\la+1}\to D_{\la}$ is bijective. We define $D^{(2)}\colon = D_{\la+1}$.\enddemo

\begin{definition}
\label{definition:Linv} Let $D$ be a monodromy $\bT$-module, let $j_0$ be the integer appearing in part (iii) of definition \ref{definition:monodmod} and $D= D^{(1)}\oplus D^{(2)}$ be the decomposition into $\p$-stable free rank-1 $\bT_K$-submodules of $D$ as in lemma \ref{lemma:phidecomp}. The $\cL$-invariant $\cL(D)$ of $D$ is defined as the unique element in $\bT_K$ such that $x-\cL(D) N(x)\in F^{j_0}D_K$ for every $x\in D^{(2)}$. 
\end{definition}

\begin{remark}
\label{remark:Lbasechange}
{\rm (a) Let $\ka'$ be a perfect extension field of $\ka$ and let $K' \supseteq K$ be the quotient field of $W(k')$. Let $D$ be a $\bT$-object in $MF_K(\p, N)$. Then $D$ is a monodromy $\bT$-module if and only if $D_L = D \otimes_K L$ is a monodromy $\bT$-module in $MF_L(\p, N)$ and that in this case the $\cL$-invariants coincide $\cL(D)=\cL(D_L)$.\\
(b) The notion of a monodromy module as introduced in (\cite{mazur}, section 9) is not quite adequate for the application to modular forms. In fact with the notation as in (\cite{mazur}, section 12) the filtered $(\p,N)$-module $D_p(f)$ considered there is a two-dimensional $F$-module but the $\p$-action is $F$-linear (not $\sigma_F$-linear) and therefore $D_p(f)$ is neither a (two-dimensional) monodromy module in $MF_F(\p, N)$ nor in $MF_{\Q_p}(\p, N)$ if $F \neq \Q$ in the sense of (\cite{mazur}, sect.\ 9).}
\end{remark}

\section{Covergent filtered $F$-isocrystals and Dieu\-donn\'e modules}
\label{sec:isocrystals}

Let $K$ be as in section \ref{sec:phin} and assume further that $K$ is of absolute ramification index $1$, i.e.\ $K = K_0$. Let $\cO_K$ be its ring of integers. Let $Z$ be a $p$-adic formal $\cO_K$-scheme i.e.\ a formal $\cO_K$-scheme locally of finite type in the sense of (\cite{bert}, (0.2.1)). We assume that $Z$ is analytically smooth (\cite{ogus}, p.772), so that the  rigid analytic $K$-variety $Z^{an}$ associated to $Z$ is smooth (for the construction of $Z^{an}$ we refer to \cite{bert}, (0.2.1)).

We want to define the notion of a filtered $F$-isocrystal on $Z$. It is a convergent $F$-isocrystal $\cE$ together with a filtration on the coherent $\cO_{Z^{an}}$-module $E^{an}$ on $Z^{an}$ associated to $\cE$ satisfying ``Griffith transversality'' with respect to the connection. To give a precise definition we first recall briefly the notion of an $F$-isocrystal from (\cite{ogus}, 2.17).

An enlargement of $Z$ is a pair $(T,z_T)$ consisting of a flat formal $\cO_K$-scheme $T$ and a morphism of formal $\cO_K$-schemes $z_T:T_0\lra Z$, where
the subscript $0$ means the reduced closed subscheme of the closed subscheme 
of $T$ defined by the ideal $p\cO_T$ (\cite{ogus}, 2.1). 

\begin{definition}
\label{def:isocrystals}
A convergent isocrystal $\cE$ on $Z$ consists of the following data:\\
(a) For every enlargement $T = (T,z_T)$ of $Z$ a coherent $\cO_T\otimes_{\cO_K}K$-module $\cE_T$.\\
(b)  For every morphism of enlargements $g:(T',z_{T'})\lra (T,z_T)$ an 
isomorphism of $\cO_{T'}\otimes_{\cO_K}K$-modules
\[
\theta _g: g^*(\cE_T)\lra \cE_{T'}.
\]
The collection of isomorphisms $\{\theta _g\}$ is required to satisfy 
the cocycle condition.
\end{definition}

If $T$ is an enlargement of $Z$ then $\cE_T$ may be interpreted 
as a coherent sheaf $E^{an}_T$ on the rigid space 
$T^{an}$ (\cite{ogus}, 1.5). If we assume that $T$ is analytically smooth over $\cO_K$ then there is a natural integrable connection 
\begin{equation}
\label{eqn:conn}
\nabla_T:E^{an}_T\lra E^{an}_T\otimes \Omega^1_{T^{an}}
\end{equation}
(cf.\ \cite{ogus}, 1.20, 2.81).

Note that the notion of an isocrystal on $Z$ depends only on the $\ka$-scheme $Z_0$. More precisely there is an equivalence of categories between convergent 
isocrystals on $Z$ and on $Z_0$. We identify an isocrystal on $Z$ and $Z_0$ in the following. Let $F = F_{Z_0}$ denote the absolute Frobenius of $Z_0$. 

\begin{definition}
\label{def:F_isocrystals}
A convergent $F$-isocrystal on $Z$ is a convergent isocrystal $\cE$ on $Z$ together with an isomorphism of isocrystals $\Phi:F^*\cE\to \cE$.
\end{definition}

We define filtered $F$-isocrystals only in the case where $Z$ is analytically smooth over $\cO_K$ (which we assume from now on) so that $E^{an} \colon = E^{an}_Z$ carries a natural connection $\nabla$.

\begin{definition}
\label{def:fiso}
A filtered convergent $F$-isocrystal on $Z$ consists of an $F$-isocrystal $\cE$ together with an exhaustive and separated decreasing filtration $F^{\bullet} E^{an}$ of coherent $\cO_{Z^{an}}$-submodules such that $\nabla(\cF^i) \seq \cF^{i-1}\otimes_{\cO_Z^{an}} \Omega^1_{Z^{an}}$ for all $i$.
\end{definition}

The category of convergent filtered isocrystals on $Z$ is an additive tensor category.

\begin{exam}
\label{exam:fic}
{\rm (a) The category of convergent filtered isocrystals on $Z$ is an additive rigid tensor category. The assignment $T \mapsto \cO_T\otimes K$ -- denoted by $\cO_Z$ -- with the canonical Frobenius and filtration given by $F^i = \cO_{Z^{an}}$ for $i\leq 0$ and $F^i = 0$ for $i>0$ is the identity object.\\
(b) A filtered Frobenius module on $Z =\spf(\cO_K)$ is a filtered Frobenius ($K$-)module in the sense of section \ref{sec:phin}.\\
(c) If $f: X \to Z$ is a smooth proper morphism of $p$-adic formal schemes then the $F$-isocrystal $R^q f_* \cO_{X/K}$ defined in (\cite{ogus}, 3.1) -- built using crystalline cohomology sheaves $\otimes K$ -- is a convergent filtered $F$-isocrystal in a natural way. In fact the associated coherent $\cO_{Z^{an}}$-module $(R^q f_* \cO_{X/K})^{an}$ is isomorphic to the relative de Rham cohomology $\bH^q_{DR}(X^{an}/Z^{an}) = R^q f_* \Omega^{\bullet}_{X^{an}/Z^{an}}$ and the connection (\ref{eqn:conn}) coincides with the Gauss-Manin connection (\cite{ogus}, 3.10). The filtration on $(R^q f_* \cO_{X/K})^{an} \cong \bH^q_{DR}(X^{an}/Z^{an})$ is the Hodge filtration.}
\end{exam}

Let $\fX \to \spec(\cO_K)$ be a proper semistable curve with connected fibers. We assume that the generic fiber $X$ is smooth and projective and that the irreducible components $C_1, \ldots, C_r$ of the special fiber $C$ are all smooth and geometrically connected and that there is more than one of them. We assume moreover that the singular points of $C$ are all $\ka$-rational ordinary double points. Let $\cE$ be a convergent filtered $F$-isocrystal on the formal completion $\wh{\fX}$ of $\fX$. The rigid analytic space $X^{an}$ associated to $X$ (see \cite{bert}, 0.3) can be identified with $\wh{\fX}^{an}$ and the coherent $\cO_{X^{an}}$-module $E^{an}$ with connection and filtration is defined algebraically i.e.\ there exists a coherent locally free $\cO_X$-module $E$ with connection $\nabla:E\lra E\otimes_{\cO_X} \Omega^1_X$ and filtration $F^{\bullet} E$ by $\cO_X$-submodules such that $\nabla(\cF^i)\seq \cF^{i-1}\otimes_{\cO_X} \Omega^1_X$ and the data yield the connection and filtration after passing to $E^{an}$.

In \cite{cole}, R.\ Coleman has defined a structure of a filtered $(\p,N)$-module on the de Rham cohomology group $H^1_{DR}(X, E)$ by using $p$-adic integration. We briefly recall his construction (for details we refer to \cite{cole} and \cite{coleman_iovita}).

The Hodge filtration on $H^1_{DR}(X, E)$ is defined as 
\begin{equation}
\label{eqn:hodge}
F^i H^1_{DR}(X, E) = \Image(H^1(X, F^i E \stackrel{\nabla}{\to} F^{i-1} E\otimes \Omega^1_X) \lra 
H^1(X, E \otimes \Omega^{\bullet})).
\end{equation}
The Frobenius and monodromy operator are defined using the admissible covering of $X^{an}$ given by the tubes around components of $C$. For that let $\cG$ be the (oriented) intersection graph of $C$. The vertices $\vertices(\cG)$ are the irreducible components of $C$ and the oriented edges $\edges(\cG)$ are triples $(x, C_i, C_j)$ where $x$ is a singular point of $C$ and $C_i, C_j$ are the two components on which $x$ lies. If $e = (x, C_i, C_j)$ is an edge we set $o(e) = C_i, t(e) = C_j$ and let $\bar{e}$ be the opposite edge $(x, C_j, C_i)$. For a vertex $v = C_i$ of $\cG$ we let $U_v = \red^{-1}(C_i)$ be the tube associated to it where $\red: X^{an} \to C(\kbar)$ is the reduction map. Similarly for an edge $e = (x, C_i, C_j)$ we let $A_e$ be the wide open annulus $\red^{-1}(x) = U_{o(e)} \cap U_{t(e)}$. The orientation of $e$ induces an orientation on $A_e$. The open cover $\{U_v\}_{v\in\vertices(\cG)}$ of $X^{an}$ is admissible. The corresponding Mayer-Vietoris sequence yields a short exact sequence
\begin{equation}
\label{eqn:mayer_vietoris}
\begin{array}{c}
0\lra (\bigoplus _{e\in \edges(\cG)}\,H^0_{DR}(A_e, E^{an}))^{-}/
(\bigoplus _{v\in \vertices(\cG)}\, H^0_{DR}(U_v, E^{an}))
\stackrel{\iota}{\lra} H^1_{DR}(X, E)\\
\stackrel{\gamma}{\lra}\Ker(\bigoplus _{v\in \vertices(\cG)}\,H^1_{DR}(U_v,E^{an}) \to \bigoplus_{e\in \edges(\cG)}\,H^1_{DR}(A_e, E^{an}))\lra 0.
\end{array}
\end{equation}
The superscript $-$ indicates the subspace of $\bigoplus _{e\in \edges(\cG)}\,H^0_{DR}(A_e, E^{an})$ consisting of elements $\{f_e\}_{e\in \edges(\cG)}$ with $f_{\bar{e}} = - f_e$ for all $e\in \edges(\cG)$.

The left and right terms in the exact sequence (\ref{eqn:mayer_vietoris}) have natural Frobenii induced from the Frobenius on $\cE$. Moreover the map $\iota$ admits a natural left inverse $s$ defined as follows: by identifying $H^1_{DR}(X, E)$ with the Cech hypercohomology of the covering $\{U_v\}_{v\in\vertices(\cG)}$ we can represent elements $\omega\in H^1_{DR}(X, E)$ as pairs of collections $(\{\omega _v\}_{v\in v(G)}, \{f_e\}_{e\in e(G)})$ where $\omega _v\in (E^{an}\otimes\Omega ^1_{X^{an}})(U_v)$ and $f_e\in E^{an}(A_e)$ are such that $f_{\bar{e}}=-f_e$ and $\omega _{o(e)}|_{A_e}-\omega_{t(e)}|_{A_e}=\nabla(f_e)$ for all $e\in \edges(\cG)$. Then $s(\omega)$ will be represented by the family $\{g_e\}_{e\in \edges(\cG)}$ where $g_e = f_e-(\lambda_{o(e)}|_{A_e}-\lambda_{t(e)}|_{A_e})$. Here, if $v\in \vertices(\cG)$ then $\lambda _v$ is a $p$-adic integral of $\omega _v$. It is well defined up to a rigid horizontal section of $E^{an}|_{U_v}$ once we have fixed a branch of the $p$-adic logarithm (which we assume from now on). The map $s$ defines a splitting of the sequence (\ref{eqn:mayer_vietoris}) and the Frobenius $\Phi$ on $H^1_{DR}(X, E)$ is defined by requiring it to be compatible with the Frobenii on the right and left term and the splitting.

The monodromy operator is defined using residues. There is a natural residue map
\[
\Res_e : H^1_{DR}(A_e,E^{an})\lra H^0_{DR}(A_e,E^{an}) \cong (E^{an}|_{A_e})^{\nabla = 0}.
\]
(to define it one has to assume that $E^{an}|_{A_e}$ has a basis of horizontal sections; however this is always the case as it is shown in \cite{coleman_iovita}, section 3.2). The monodromy operator is defined as the composition
\begin{equation}
\label{eqn:monr}
\begin{array}{c}
N =\iota\circ (\bigoplus_e\, \Res_e) : H^1_{DR}(X, E)\lra H^1_{DR}(X, E)\\
\omega \mapsto \iota(\{\Res_e(\omega|_{A_e})\}_{e\in\edges(\cG)})
\end{array}
\end{equation}
It satisfies the relation $N \Phi = p\Phi N$.

We also need to consider the de Rham cohomology of an open subscheme $U$ of $X$ which is the complement of a finite number of points which specialise to smooth points on $C$. More precisely let $U = X-S$ where $S$ is a finite set of $K$-rational points of $X$ which -- considered as points on $\fX$ -- are all smooth and which specialise to pairwise different (smooth) points on $C$. Then one can define in a similar way a structure of a filtered $(\p,N)$-module on $H^1_{DR}(U, E)$. In fact if we replace $X$ by $U$ and $U_v$ by $U_v-S$ in the sequence (\ref{eqn:mayer_vietoris}) and in (\ref{eqn:monr}) the sequence (\ref{eqn:mayer_vietoris}) is also exact and the monodromy operator is defined again by (\ref{eqn:monr}). Moreover $\iota$ has a left inverse defined as before. Note that the term on the left in (\ref{eqn:mayer_vietoris}) remains the same since $H^0_{DR}(U_v, E^{an})$ does not change if we remove a finite number of points from $U_v$. However to define the Frobenius on the right term (and thus on $H^1_{DR}(U, E)$) is more involved as one has to work with logarithmic isocrystals. Associated to the divisor $\sum_{P\in S} \, P$ on $\fX$ is a fine logarithmic structure $M_S$ on $\fX$ and $\wh{\fX}$ (see \cite{kato}). By pulling back $\cE$ under the canonical morphism of formal log-schemes 
$j:(\wh{\fX}, M) \to (\wh{\fX},\, \mbox{trivial log-structure})$ one gets a convergent log-$F$-isocrystal (a log-isocrystal is defined in terms of log-enlargements -- a log-version of enlargements). Then the de Rham cohomology $H^1_{DR}(U_v-S, E^{an})$ can be described in terms of log-crystalline cohomology with coefficients in $j^*\cE$ of the component $C_i = v$ of $C$ and thus carries a Frobenius. For details we refer to \cite{coleman_iovita}.
 
The Gysin sequence
\begin{equation}
\label{eqn:gysin}
\begin{CD}
0 @>>> H^1_{DR}(X, E) @>>> H^1_{DR}(U, E) @> \oplus_{x\in S}\,\Res_x >> \bigoplus_{x\in S} \, {\cE}_x[1]
\end{CD}
\end{equation}
becomes then an exact sequence of filtered $(\p, N)$-modules. Here ${\cE}_x[1]$ is the stalk of $\cE$ at $x$; by Example \ref{exam:fic}(b) it is a filtered Frobenius module.

In \cite{coleman_iovita} the following technical assumption on a convergent
$F$-isocrystal $\cE$ on $C$ is made in order to compare Coleman's definition of the $(\p,N)$-module $H^1_{\DR}(X, \cE)$ with Faltings' log-crystalline cohomology. Let $i:x\rightarrow C$ be the embedding of a closed point in $C$ and $j:C_v^0\rightarrow C$ be the embedding of the complement of the singular points of $C$ in a component. Then 
\begin{assumption}
\label{ass:frob}
The isocrystal $H^0_{crys}(x, i^*(\cE))$ is isotypical and its slope does not occur in the set of slopes of $H^1_{crys}(C_v^0, j^*(\cE))$.
\end{assumption}

Consider now the case where $\cE$ is a convergent filtered isocrystal of the type considered in Example \ref{exam:fic}(c). More precisely let $f:\fY \to \fX$ be a smooth proper morphism, $Y = \fY_{K}$ the generic fiber of $\fY$ and let $\bH^q_{DR}(Y/X) \colon = R^q f_* \cO_{\wh{\fY}/K}$ be the convergent filtered $F$-isocrystal of \ref{exam:fic}(c). The results of \cite{faltings1, faltings2} and \cite{coleman_iovita} imply the following.

\begin{theorem}
\label{theorem:comp} 
(a) The representation $H^1_{\et}(\Xbar,R^q f_* \Q_p))$ is semistable. If assumption \ref{ass:frob} is satisfied for the isocrystal $\bH^q_{DR}(Y/X)$ then $D_{\st}(H^1_{\et}(\Xbar,R^q f_* \Q_p))$ is canonically isomorphic as a filtered $(\p, N)$-module to $H^1_{DR}(X,\bH^q_{DR}(Y/X))$.\\
(b) More generally let $S$ be a finite set of smooth sections of $f:\fX \to \spec(\cO_K)$ which specialise to pairwise different (smooth) points on $C$ and let $U = X -S$, $\Ubar = U\otimes_{K} \Kbar$ and $Y_{\bar{x}}$ be the geometric fiber of $f:Y \to X$ over $x\in S$. Then we have an exact sequence of semistable Galois representations
\[
0 \lra H^1_{\et}(\Xbar,R^q f_* \Q_p))\lra H^1_{\et}(\Ubar,R^q f_* \Q_p))\lra \bigoplus_{x\in S} H^q_{\et}(Y_{\bar{x}}, \Q_p(-1))
\]
which is -- after applying the functor $D_{\st}$ -- isomorphic to the sequence {\rm (\ref{eqn:gysin})} (for $\cE = \bH^q_{DR}(Y/X)$.
\end{theorem}

\section{Convergent filtered $F$-isocrystals associated to representations of $\GL_2$}
\label{sec:gltwo}

Let $M_2$ be the algebra scheme of $2\times 2$-matrices and $\GL_2$ the group scheme of invertible elements in $M_2$. By $\cH_p$ we denote the $p$-adic upper half plane over $\Q_p$; it is a rigid analytic $\Q_p$-variety whose $\C_p$-valued points are $\cH_p(\C_p) = \C_p -\Q_p$. We denote by $\wh{\cH}$ the canonical formal model of $\cH_p$ over $\Zur$. We have a left $\GL_2(\Q_p)$-action through linear transformations on $\cH_p$. The set $\cH_p(\Q_{p^2})$ will be often identified with the set of $\Q_p$-algebra homomorphisms $\Hom(\Q_{p^2}, M_2(\Q_p))$ as follows. Any $\psi\in \Hom(\Q_{p^2}, M_2(\Q_p))$ defines an action of $\Q_{p^2}^*$ on $\cH_p(\Q_{p^2})$ and the point $z\in \Q_{p^2}-\Q_p$ corresponding to $\psi$ is characterized by the property
\begin{equation}
\label{eqn:puh}
\psi(a) \left(\begin{array}{c} z \\ 1\end{array}\right) = a \left(\begin{array}{c} z \\ 1\end{array}\right)\qquad\forall\, a\in \Q_{p^2}.
\end{equation}

For an algebraic group $\GG$ over $\Q_p$ we let $\Rep_{\Q_p}(\GG)$ be the category of finite-dimensional $\Q_p$-rational representations $\rho: \GG \to \GL(V)$. In this section we construct for every representation $V$ in $\Rep_{\Q_p}(\GL_2\times \GL_2)$ a filtered $F$-isocrystal $\cE(V)$ on $\wh{\cH}$ and for every $\Q_{p^2}$-rational point $\psi\in \Hom(\Q_{p^2}, M_2(\Q_p))$ a filtered Frobenius module $V_{\Psi}\in MF_{\Qur}(\p, N)$ which turns out to be the fiber of $\cE(V)$ at $\Psi$ (compare also \cite{rap-zink}, 1.31 for a related construction).

For an object $(V,\rho: \GL_2\times \GL_2\to \GL(V))$ in $\Rep_{\Q_p}(\GL_2\times \GL_2)$ we also use the notation $(V, \rho_1, \rho_2)$ where $\rho_1: \GL_2 \to \GL(V)$ (resp.\ $\rho_2$) is the restriction of $\rho$ to the first (resp.\ second) factor. If $(V_1,\rho_1), (V_2,\rho_2)$ are two representations of $\GL_2$ we define $V_1 \odot V_2 =  (V_1,\rho_1)\odot (V_2,\rho_2)$ to be the $\GL_2\times \GL_2$-representation
\[
V_1 \odot V_2 = (V_1\otimes V_2, \rho_1\otimes 1_{V_2},1_{V_1} \otimes \rho_2).
\]
For $V\in \Rep_{\Q_p}(\GL_2)$ and $m\in \Z$ we write $V\{m\} \colon = V \odot {\det}^{\otimes^m}, \{m\}V \colon = {\det}^{\otimes^m}\odot V$.

A representation $(V, \rho)\in \Rep_{\Q_p}(\GL_2)$ is said to be pure of weight $n$ if every element $a$ of the center $Z(\GL_2) = \BGG_m$ acts by multiplication with $a^n$. In general $(V, \rho)$ can be written as
\begin{equation}
\label{eqn:weight}
(V, \rho) = \bigoplus_{n\in \Z}\, (V^{(n)}, \rho^{(n)})
\end{equation}
where $(V^{(n)}, \rho^{(n)})$ is pure of weight $n$.

For $n\geq 0$ let $\cP_n$ denote the vector space of polynomials of degree $\leq n$ over $\Q_p$ with right $\GL_2$-action
\[
P(X) \cdot A  = (cX+d)^n P(\frac{aX+b}{cX+d}) \qquad A = \mat{a}{b}{c}{d}\in \GL_2, P(X)\in \cP_n. 
\]
Therefore the dual $V_n \colon = \cP_n^{\vee} = \Hom_{\Q_p}(\cP_n, \Q_p)$ has a left $\GL_2$-action given by 
\begin{equation}
\label{eqn:vn}
(A \cdot R)(P(X)) = R(P(X) \cdot A) \qquad \forall\, P(X)\in \cP_n.
\end{equation}
We need the following simple fact (see \cite{deligne2}, Proposition 3.1).

\begin{lemma}
\label{lemma:gv1}
Any object of $\Rep_{\Q_p}(\GL_2)$ is a direct summand of a sum of representations of the form $V_1^{\otimes^m}\otimes (V_1^{\vee})^{\otimes^n}, m, n\geq 0$.
\end{lemma}

For $(V, \rho_1, \rho_2)$ in $\Rep_{\Q_p}(\GL_2\times \GL_2)$ let $\phi_V = \rho_2(\mat{0}{p}{1}{0}): V \to V$. Note that if $(V, \rho_1, \rho_2)$ is of the form $V\{m\}$ with $V\in Ob(\Rep_{\Q_p}(\GL_2))$ then $\phi_{V\{m\}} = p^m 1_V$.

For $\Psi\in \Hom(\Q_{p^2}, M_2(\Q_p))$ the underlying vector space of the filtered Frobenius module $V_{\Psi}$ is $V_{\Qur}$. The filtration is defined in terms of the first $\GL_2$-action $\rho_1$ solely. Consider first the case where $(V, \rho_1)$ is pure of weight $n$. Let $V_{\Qur,i}$ be the subspace of $v\in V_{\Qur}$ such that $\rho_1(\Psi(a))(v) = a^i \sigma(a)^{n-i} v$ for all $a\in \Q_{p^2}$ 
and
\[
F^i_{\Psi} V_{\Qur} \colon = \bigoplus_{j\geq i} V_{\Qur,j}.
\]
For arbitrary $V$ we use the decomposition (\ref{eqn:weight}) and set
\[
F^i_{\Psi} V_{\Qur} \colon = \bigoplus_{n\in \Z}\,F^i_{\Psi} V^{(n)}_{\Qur}
\]
Then $F^{\bullet}_{\Psi}$ is an exhaustive and separated filtration on $V_{\Qur}$ and we define $V_{\Psi}$ to be $V_{\Qur}$ together with Frobenius $\Phi = \phi_V \otimes \sigma$ and filtration $F^{\bullet}_{\Psi}$. The functor $V \mapsto V_{\Psi}$ from $\Rep_{\Q_p}(\GL_2\times \GL_2)$ to $MF_{\Qur}(\p, N)$ is an exact tensor functor.

The convergent filtered $F$-isocrystal $\cE(V)$ associated to $(V, \rho_1, \rho_2)$ in $\Rep_{\Q_p} \break ( \GL_2 \times \GL_2)$ is defined as follows. As an isocrystal it is just $V\otimes_{\Q_p} \cO_{\wh{\cH}}$ where $\cO_{\wh{\cH}}$ is the isocrystal of example (\ref{exam:fic}) (a). The Frobenius is given by $\phi_V \otimes \Phi_{\cO_{\cH_p}}$. The filtration is again defined in terms of $\rho_1: \GL_2 \to \GL(V)$ only. Consider first the case $(V, \rho_1) = V_1$. We have a canonical map of sheaves on $\cH_p$
\begin{equation}
\label{eqn:ev}
V_1 \otimes \cO_{\cH_p} \lra \cO_{\cH_p}, R \otimes f \mapsto (z\mapsto R(X-z)\, f(z)),
\end{equation}
(where $R(X-z) \colon = R(X) - z R(1)$). The filtration of $V_1\otimes_{\Q_p} \cO_{\wh{\cH}}$ is given by $F^0 = V_1\otimes_{\Q_p} \cO_{\wh{\cH}}$, $F^2 = 0$ and $F^1$ is the kernel of (\ref{eqn:ev}). We get an induced filtration on $V_1^{\otimes^m}\otimes (V_1^{\vee})^{\otimes^n})\otimes_{\Q_p} \cO_{\wh{\cH}}$ for all $m, n\geq 0$ and thus by Lemma \ref{lemma:gv1} on $\cE(V)^{an}$ for an arbitrary representation $V$ in $\Rep_{\Q_p}(\GL_2\times \GL_2)$. 

Again the assignment $V \mapsto \cE(V)$ is an exact tensor functor. For an even integer $n =2m$ we write $\cV_n$ for $\cE(V_n\{m\})$.

\begin{lemma}
\label{lemma:stalk} 
For $(V, \rho_1, \rho_2)$ in $\Rep_{\Q_p}(\GL_2\times \GL_2)$ and $\Psi\in \Hom(\Q_{p^2}, M_2(\Q_p))$ we have $\cE(V)_{\Psi} \cong V_{\Psi}$, where $\cE(V)_{\Psi}$ is the stalk of $\cE(V)$ at $\Psi$.
\end{lemma}

{\em Proof.} That the Frobenii are the same is obvious. To see that the filtrations also agree it is enough, by Lemma \ref{lemma:gv1}, to consider the case $(V, \rho_1) = V_1$. Note that $(V_1\otimes_{\Q_p} \cO_{\wh{\cH}})_{\Psi}$ is generated by $ev_z\in (V_1)_{\Qur} = \Hom(\cP_1, \Qur), \break ev_z(P(X)) \colon= P(z)$ and that 
\[
(\Psi(a)\cdot ev_z)(P(X)) = a P(z) = a (ev_z)(P(X))\qquad \forall\, a\in \Q_{p^2}^*, P(X) \in \cP_1.
\]
Here $z\in \Q_{p^2}-\Q_p$ is the point corresponding to $\Psi$ via (\ref{eqn:puh}).
\enddemo

The filtered isocrystals $\cE(V)$ constructed above also descend to isocrystals on Mumford curves. Let $\G\seq \SL_2(\Q_p)$ be a discrete cocompact subgroup  and $X= X_{\G}$ the associated Mumford curve over $\Qur$ so that $X_{\G}^{an} = \G \backslash \cH_p$. Assume first that $\G$ is torsionfree. The action of $\GL_2$ on $\cE(V)$ is compatible with the action on $\cH_p$ and therefore $\cE(V)$ descends to a filtered isocrystal on $X_{\G}$ which will be also denote by $\cE(V)$. We let $E(V)$ be the associated coherent locally free $\cO_{X_{\G}}$-module $E$ with connection and filtration. 

We will give now a concrete description of the $(\p, N)$-module structure on $H^1_{DR}(X_{\G}, E(V))$. For this we have to recall the description of the special fiber $\wh{\cH}_s$ of $\wh{\cH}$. Let $\cT$ be the Bruhat-Tits tree for $\PGL_2(\Q_p)$; the vertices are homothety classes of lattices in $\Q_p^2$. Two vertices $v_1, v_2$ are adjacent if one can choose representing lattices $M_1, M_2$ with $M_1 \subsetneq M_2 \subsetneq p M_1$. An oriented edge of $\cT$ is a pair $e = (v_1, v_2)$ where $v_1, v_2$ are adjacent vertices. We set $o(e) = v_1, t(e) = v_2$ and $\bar{e} = (v_2, v_1)$.

The set of components of $\wh{\cH}_s$ is in one-to-one correspondence to the set vertices $\cV(\cT)$ of $\cT$, each component being isomorphic to the projective line over $\F_p$. We write $\{ \PP_v^1\}_{v\in \cV(\cT)}$ for the set of components of $\wh{\cH}_s$. The singular points of $\wh{\cH}_s$ are ordinary $\F_p$-rational double points; they correspond to the set of (unoriented) edges of $\cT$; two components $\PP_{v_1}^1, \PP_{v_2}^1$ intersect if and only if $v_1, v_2$ are adjacent and in this case transversally in an $\F_p$-rational point. The curve $X_{\G}$ has a canonical semistable model $\cX_{\G}$ whose formal completion is $\G\backslash \wh{\cH}$. Hence the special fiber of $\cX_{\G}$ is $\G\backslash \wh{\cH}_s$. 

Let $C^0(V_{\Qur})$ (resp.\ $C^1(V_{\Qur})$) be the set of maps $f:\cV(\cT) \to V_{\Qur}$ (resp.\ maps $f: \edges(\cT) \to V_{\Qur}$ such that $f(\bar{e}) = -f(e)$ for all $e\in \edges(\cT)$). The short exact sequence 
\[
0 \lra V_{\Qur} \lra C^0(V_{\Qur}) \stackrel{\partial}{\lra} C^1(V_{\Qur}) \lra 0
\]
(where $\partial(f)(e) = f(o(e)) - f(t(e))$) yields an isomorphism 
\begin{equation}
\label{eqn:verb}
\ep: C^1(V_{\Qur})^{\G}/C^0(V_{\Qur})^{\G} \cong H^1(\G, V_{\Qur}).
\end{equation}
From the above description of $\cX_{\G}$ we see that the left term in the sequence (\ref{eqn:mayer_vietoris}) is equal to $C^1(V_{\Qur})^{\G}/C^0(V_{\Qur})^{\G}$ and therefore can be replaced by $H^1(\G, V_{\Qur})$. In this situation we write $P$ instead of the splitting $s$ introduced in the last section. Explicitely it is given as follows. We identify $H^1_{DR}(X_{\G}, E(V))$ with the space of $V$-valued $\G$-invariant meromorphic differential forms of the second kind on $\cH_p$ (modulo exact forms). Given such a form $\omega$ we let $F_{\omega}$ be a primitive of it (see \cite{shalit1}). Then $P$ is given by
\begin{equation}
\label{eqn:icole}
P:H^1_{DR}(X_{\G}, E(V))\lra H^1(\G, V_{\Qur}), \omega \mapsto (\ga\mapsto \ga(F_{\omega}) - F_{\omega})
\end{equation}

The monodromy operator $N$ is the composite of $\iota \circ (-\ep) \circ I$ where $I$ is Schneider's integration map (see \cite{shalit3})
\begin{equation}
\label{eqn:isch}
I: H^1_{DR}(X_{\G}, E(V)) \lra C^1(V_{\Qur})^{\G}, \omega \mapsto (e\mapsto \Res_e(\omega)).
\end{equation}
As for the Frobenius we first note that under the isomorphism (\ref{eqn:verb}) we obtain a Frobenius on $H^1(\G, V_{\Qur})$. It is the one induced by the map $\phi_V \otimes \sigma$ on the coefficients $V_{\Qur}$. 

We will assume from now on that $\G$ is arithmetic (see \cite{shalit1}). We need the following lemma to characterise the Frobenius.

\begin{lemma} Assume that $\G$ is arithmetic. Then sequence 
\[
0\lra H^1(\G, V_{\Qur}) \stackrel{\iota}{\lra} H^1_{DR}(X_{\G}, E(V))\stackrel{I}{\lra} C^1(V_{\Qur})^{\G}/C^0(V_{\Qur})^{\G} \lra 0
\]
is exact.
\end{lemma}

{\em Proof.} We may assume that $V$ is an irreducible representation of $\GL_2$. Then $V \cong V_n \otimes det^{\otimes^m}$ for some $m,n\in\Z, n\geq 0$, hence $V \cong V_n$ as a $\G$-module. The assertion follows now from (\cite{shalit1} 3.9) and (\cite{shalit3} 1.6).
\enddemo

Hence the sequence 
\begin{equation}
\label{equation:NNex}
H^1_{DR}(X_{\G}, E(V))\stackrel{N}{\lra}H^1_{DR}(X_{\G}, E(V))\stackrel{N}{\lra}H^1_{DR}(X_{\G}, E(V))
\end{equation} 
is exact and it follows that there is a unique Frobenius operator on $H^1_{DR}(X_{\G}, \break E(V))$ satisfying $N \p = p \p N$ and which is compatible (with respect to $\iota$ and $P$) with the Frobenius on $H^1(\G, V_{\Qur})$.

It is easy to see that if $(V, \rho_2)$ is pure of weight $n$ then the isocrystal $\cE(V)$ on $X_{\G}$ satisfies assumption \ref{ass:frob} (in fact with the notation there $H^0_{crys}(x, i^*(\cE(V)))$ and $H^1_{crys}(C_v^0, j^*(\cE(V)))$ are isotypical of slope $n/2$ and $n/2 + 1$ respectively).

If $\G$ is not torsionfree we define a filtered $(\p,N)$-module $H^1_{DR}(X_{\G},\break E(V))$ as follows. We choose a free normal subgroup of finite index $\G'\seq \G$. The group $\G/\G'$ acts then on $H^1_{DR}(X_{\G'}, E(V))$ as automorphisms of a filtered $(\p,N)$-module and thus the filtered $(\p,N)$-module structure on $H^1_{DR}(X_{\G'}, E(V))$ induces one on $H^1_{DR}(X_{\G}, E(V)) \colon = H^1_{DR}(X_{\G}, E(V))^{\G/\G'}$. This construction is clearly independent of the choice of $\G'$.

In the case of an open subscheme $U = X_{\G} -S$ of the type considered in the last section the description of the $(\p,N)$-module structure on $H^1_{DR}(U, E(V))$ is very similar. The monodromy operator is defined as before. Let $\pi: \cH_p \to X_{\G}^{an}$ be the canonical projection. We identify $H^1_{DR}(U, E(V))$ with the space of $V$-valued $\G$-invariant meromorphic differential forms on $\cH_p$ which are of the second kind when restricted to $\pi^{-1}(U)$. The left inverse $P: H^1_{DR}(U, E(V))\to H^1(\G, V_{\Qur})$ of $\iota$ is defined by the same formula (\ref{eqn:icole}). We choose for each $x\in S$ a point $\Psi_x\in \Hom(\Q_{p^2}, M_2(\Q_p))$ with $\pi(\Psi_x) = x$. One can easily verify that there is a unique Frobenius on $H^1_{DR}(U, E(V))$ such that the Gysin sequence
\begin{equation}
\label{eqn:gysin2}
\begin{CD}
0 @>>> H^1_{DR}(X, E(V)) @>>> H^1_{DR}(U, E(V)) @> \oplus_{x\in S}\,\Res_x >> \bigoplus_{x\in S} \, V_{\Psi_x}[1]
\end{CD}
\end{equation}
is a sequence of $(\p,N)$-modules and such that $P$ is compatible with Frobenii.

\section{Representations and Dieudonn\'e module attached to modular forms}
\label{sec:dieumod}

Let $f_{\infty}$ be a modular form of a fixed even weight $k\geq 4$ and level $N$ and assume that the prime number $p$ divides $N$ exactly once. Under certain conditions on the level one can associate to $f_{\infty}$ a modular form $f$ on a Shimura curve via the Jacquet-Langlands correspondence and the Galois representation associated to $f_{\infty}$ can be identified with a direct summand, denoted by $V_p(f)$, in the $p$-adic cohomology group of a Kuga-Sato variety over the Shimura curve. In this section we will define these representations and show -- using the results of the previous sections and the Theorem of Cerednik-Drinfeld -- that they are semistable as representations of the decomposition group at $p$ and describe explicitely their Dieudonn{\'e} modules using $p$-adic integration. 

To begin with we fix some notation. Put $n=k-2$, $m =\frac{n}{2}$. For the rest of this paper, $N^{-}$ denotes a positive squarefree integer with an odd number of prime divisors none of which equals $p$, $N^+$ denotes a positive integer relatively prime to $p N^{-}$ and $N$ the product $p N^{-} N^{+}$. Let $\cB$ be the indefinite quaternion algebra over $\Q$ of discriminant $p N^{-}$. We fix a maximal order $\cR^{\max}$ in $\cB$ and an Eichler order of level $N^{+}$ contained in $\cR^{\max}$. Associated to these data is a Shimura curve $X = X_{N^+, p N^{-}}$ whose description as a coarse moduli scheme we briefly recall from \cite{boutot_carayol}.

\begin{definition}
\label{definition:qmsurface}
Let $S$ be a $\Q$-scheme. An abelian surface with quaternionic multiplication (by $\cR^{\max}$) and level $N^+$-structure over $S$ (abelian surface with QM for short) is a triple $(A, \io, C)$ where\\
1. $A$ is an abelian scheme over $S$ of relative dimension $2$;\\
2. $\io:\cR^{\max}\to\End_S(A)$ is an inclusion defining an action of $\cR^{\max}$ on $A$;\\
3. $C$ is a subgroup scheme of $A$ which is locally isomorphic to $\Z/N^+\Z$ and is stable and locally cyclic under the action of $\cR$.
\end{definition}

\begin{definition}
\label{definition:cmoduli}
The Shimura curve $X = X_{N^+, pN^-}/\Q$ is the coarse moduli scheme of the moduli problem 
\[
S \mapsto \{\mbox{\,isomorphism classes of abelian surfaces with QM over $S$\,\,}\}.
\] 
\end{definition}

The scheme $X$ is a smooth projective geometrically connected curve over $\Q$. We recall its description as a Mumford curve over $\Q_{p^2}$. Let $B/\Q$ be the definite quaternion algebra over $\Q$ of discriminant $N^{-}$ and let $R$ be an Eichler $\Z[\frac{1}{p}]$-order of level $N^+$ in $B$. By fixing an isomorphism $B_p\cong M_2(\Q_p)$ the group $\G$ of elements of reduced norm 1 in $R$ can be viewed as a discrete cocompact subgroup of $\SL_2(\Q_p)$.

\begin{theorem}
\label{theorem:cer_drinf}
{\rm (\v Cerednik-Drinfeld; see \cite{boutot_carayol}, Chapitre III, 5.3.1, \cite{cerednik}, \cite{drinfeld})} We have $X_{\Q_{p^2}} \cong X_{\G}$.
\end{theorem}

We recall now the Jacquet-Langlands correspondence between cusp forms of weight $k$ and level $N$ which are $p N^-$-new and modular forms of weight $k$ on $X$. Let us first clearify the latter notion.

\begin{definition}
\label{definition:modform} 
Let $K$ be a field of characteristic $0$. A $K$-valued modular form of weight $k$ on $X$ is a global section of $\Omega_{X_K/K}^{\otimes^{m+1}}$. We denote the space of these modular forms by $M_k(X, K)$. 
\end{definition}

We remark that over the field $K = \Qur$ (or any complete field $K\seq \C_p$ containing $\Q_{p^2}$) one can give a more concrete description of $M_k(X, K)$ due to theorem \ref{theorem:cer_drinf}. A $p$-adic modular form of weight $k$ for $\G$ is a rigid analytic function $f$ on $\cH_p$ defined over $K$ such that 
\[
f(\ga z) = (cz+d)^k f(z) \qquad \mbox{for all $\ga =\mat{a}{b}{c}{d} \in \G$}.
\]
The space of these $p$-adic modular forms will be denoted by $M_k(\G) = M_k(\G,K)$. By identifying $M_k(X, \Qur)$ with the space of global section of $\Omega_{\G\backslash \cH_p}^{\otimes^{m+1}}$ we get an isomorphism
\begin{equation}
\label{eqn:rigidmodform}
M_k(\G) \lra M_k(X, \Qur), f \mapsto f(z) dz^{\otimes^{m+1}}
\end{equation}

Returning to the case of arbitrary coefficients $K$ the Jacquet-Langlands correspondence is an isomorphism (uniquely determined up to scaling)
\begin{equation}
\label{eqn:jacquet}
\begin{array}{ccc}
M_k(X, K) & \stackrel{\cong}{\longleftrightarrow} & S_k(\G_0(N), K)^{p N^{-}-new}\\
f & \stackrel{\cong}{\longleftrightarrow} & f_{\infty}
\end{array}
\end{equation} 
which is compatible with the action of the Hecke operators and Atkin-Lehner involutions.

Let $f_{\infty}\in S_k(\G_0(N), \Q)$ be a newform and $f\in M_k(X,\Q)$ the corresponding newform on $X$ and let $F = F_f$ be the finite extension of $\Q$ generated by the eigenvalues of the Hecke operators acting on $f$ (or $f_{\infty}$). Associated to $f_{\infty}$ is a two-dimensional (as $F_p$-module) $G_{\Q}$-representation $V_p(f_{\infty})$. It is defined as a direct summand of the $(n+1)$-th $p$-adic cohomology of a suitable compactification of the $n$-fold fibre product of the universal elliptic curve (with full level $N$-structure) over the modular curve $X(N)$ (see \cite{deligne1}). 

A similar construction for $f$ can be done by using the universal abelian surface over the Shimura curve. For this we have to work with a more refined moduli problem then in definition \ref{definition:cmoduli}.

\begin{definition}
\label{definition:qmlevelM} 
Let $M\geq 3$ be an integer relatively prime to $N$ and $S$ be a $\Q$-scheme.
An abelian surface with quaternionic multiplication (by $\cR^{\max}$), level $N^+$-structure and full level $M$-structure is a quadruple $(A, \io, C, \bar{\nu})$ where $(A, \io, C)$ is a triple as in \ref{definition:qmsurface} and $\bar{\nu}: (\cR^{max}/M\cR^{max})_S\to A[M]$ is a $\cR^{\max}$-equivariant isomorphism from the constant group scheme $(\cR^{max}/M\cR^{max})_S$ to the group scheme of $M$-division points of $A$.
\end{definition}

The corresponding moduli problem admits a fine moduli scheme which we denote by $X_M$ (see \cite{boutot_carayol}). Again it is a smooth projective curve over $\Q$ (however it is not geometrically connected). We have a Galois covering $q: X_M \to X$ with Galois group $\cong G/\{\pm 1\}$ where 
\[
G = G_M \colon = (\cR^{max}/M\cR^{max})^* \cong (\cR/M\cR)^* \cong \GL_2(\Z/M\Z) \cong (R/M R)^*
\]
obtained by forgetting the level $M$-structure. Let $\pi:\cA \to X_M$ be the universal abelian surface over $X_M$. The action of $\cR^{max}$ on $\cA$ induces an action of $\cB^*$ on $R^q f_* \Q_p$. We let
\[
\bL_2 \colon = \bigcap_{b\in \cB} \Ker(b- N(b): R^2 \pi_*\Q_p \to R^2 \pi_*\Q_p).
\]
It is a 3-dimensional $p$-adic local system on $X_M$. Let 
\begin{equation}
\label{eqn:lap1}
\De_m: \sym^m \bL_2 \lra (\sym^{m-2} \bL_2)(-2),
\end{equation}
be the Laplace operator associated to the non-degenerated pairing
\begin{equation}
\label{eqn:s4pair1}
(\,\,\,,\,\,\,):\bL_2\otimes\bL_2 \hookrightarrow R^2 \pi_*\Q_p\otimes R^2 \pi_*\Q_p \stackrel{\cup}{\lra} R^4 \pi_*\Q_p\stackrel{tr}{\lra} \Q_p(-2).
\end{equation}
Symbolically (\ref{eqn:lap1}) is given by
\[ 
\De_m(x_1\cdots x_m) = \sum_{1\leq i<j\leq m}\,(x_i, x_j)\, x_1\cdots \widehat{x_i}\cdots \widehat{x_j}\cdots x_m.
\]
For $n> 2$ we define $\bL_n$ as the kernel of $\De_m$. Note that we have an action of $G$ on $\bL_n$ compatible with the action on $X_M$.

\begin{definition} 
\label{definition:repsmod}
The $p$-adic $G_{\Q}$-representation attached to the space $M_k(X,\Q)$ is defined to be the representation
\[
H_p(\cM_n) = H^1(\overline{X}_M, \bL_n)^G
\]
where $\overline{X}_M = X_M\otimes_{\Q} \Qbar$.
\end{definition}

To justify the notation $H_p(\cM_n)$ note that it is the $(n+1)$-th $p$-adic realization of the motive $\cM_n$ constructed in appendix \ref{subsection:a1} (see \ref{lemma:a1hpm}). As in \cite{deligne1} or (\cite{scholl}, section 4) one can define Hecke operators and Atkin-Lehner involutions acting on $H_p(\cM_n)$ (and on $\cM_n$). Let $\TT$ be the subalgebra of $\End(S_k(\G_0(N), \Q)^{p N^{-}-new}) \cong \End(M_k(X, \Q))$ generated by the Hecke operators $T_{\ell}$ for $\ell\nmid N$ and $U_{\ell}$ for $\ell\mid N$. Then $M_k(X, \Q)$ is a free $\TT$-module of rank one. By multiplicity one to every newform $f\in M_k(X, \Q)$ corresponds a primitive idempotent $e_f\in \TT$ such that $e_f\TT \cong F_f$ and $e_f\cdot M_k(X, \Q) = F_f\cdot f$. The Betti realization $H_B(\cM_n)$ has a Hodge structure of type $(0,n+1), (n+1,0)$, $H_B(\cM_n)\otimes \C = H^{(0,n+1)} \oplus H^{(n+1,0)}$ with $H^{(0,n+1)} \cong M_k(X, \C)$. The comparison isomorphism $H_B(\cM_n)\otimes \C \cong H_p(\cM_n) \otimes_{\Q_p,\tau} \C$ -- where $\tau: \Q_p \hookrightarrow \C$ is any embedding -- shows that $H_p(\cM_n)$ is a free $\TT_{\Q_p}$-module of rank 2. 

\begin{definition} 
\label{definition:repf} 
The $p$-adic $G_{\Q}$-representation $V_p(f)$ attached to $f$ is given by $V_p(f) = e_f\, H_p(\cM_n)$.
\end{definition}

\begin{lemma}
\label{lemma:jlrep}
We have $V_p(f)\cong V_p(f_{\infty})$ as $G_{\Q}$-representations and $\TT$-modules.
\end{lemma}

{\em Proof}: The proof is similar to the one in the weight two case given by Ribet (\cite{ribet1}, lemme to th\'eor\`eme 2). The Eichler-Shimura relations imply that the traces of the Frobenii at primes $\ell\nmid N$ operating on $V_p(\p)$ and $V_p(\p')$ are given by the action of the $T_{\ell}$ and thus are the same. Cebotarev's density theorem and the theorem of Brauer-Nesbitt allow us to deduce that the semisimplifications of the $G_{\Q}$-representations $V_p(f), V_p(f_{\infty})$ are isomorphic. But in \cite{ribet2} it is shown that $V_p(f_{\infty})$ is already a simple $G_{\Q}$-representation. Thus $V_p(f)\cong V_p(f_{\infty})$.\enddemo

The main result of this section is

\begin{theorem}
\label{theorem:dieumodform}
Considered as a $G_{\Q_p}$-representation $H_p(\cM_n)$ is semistable. We have an isomorphism of $(\p,N)$-modules (canonical up to scaling)
\[
D_{\st, \Qur}(H_p(\cM_n)) = H^1_{\DR}(X_{\G}, \cV_n).
\]
\end{theorem}

{\em Proof.} For the first statement it is enough to see that $H_p(\cM_n)$ is semistable as a $G_{\Qur}$-representation. 

Over $\Qur$ the curve $X_M$ admits a $p$-adic uniformization similar to \ref{theorem:cer_drinf}. We have (see \cite{boutot_carayol}, Chapitre III, 5.3.1)
\begin{equation}
\label{eqn:puniv}
(X_M)_{\Qur}^{an} \cong \G \backslash (\cH_p \times (R/M R)^*)
\end{equation}
where $\G$ acts diagonally on $\cH_p \times (\cR/M\cR)^*$ (on the factor $(\cR/M\cR)^*$ an element $\ga\in \G$ acts by left multiplication with $\ga \mod M$). Since the orbits of the $\G$-action on $(\cR/M\cR)^*$ are the fibers of the reduced norm $Nrd: (\cR/M\cR)^* \to (\Z/M\Z)^*$ we can write (\ref{eqn:puniv}) also as 
\begin{equation}
\label{eqn:puniv2}
(X_M)_{\Qur}^{an} \cong \coprod_{(\Z/M\Z)^*} \, \G_M\backslash \cH_p
\end{equation}
where $\G_M \colon = \{\ga\in \G\mid \, ga \equiv 1 \mod M\}$, i.e.\ $(X_M)_{\Qur}$ is the disjoint union of Mumford curves $X_{\G_M}$ (indexed by the set $(\Z/M\Z)^*$). The representation (\ref{eqn:puniv}) has the advantage that the action of $G(X_M/X) \cong (\cR/M\cR)^*$ on the left hand side is more transparent (namely for $g\in (\cR/M\cR)^*$ it is given by right multiplication with $g^{-1}$ on the second factor).

To deduce \ref{theorem:dieumodform} from theorem \ref{theorem:comp} we need an explicit description of the filtered isocrystal $\bH^1_{DR}(\cA/X_M)$ (over $\Qur$). Note that $\bH^1_{DR}(\cA/X_M)$ carries a natural $\cB_{\Qur}$-action induced by the $\cR^{max}$-action on $\cA$.

On $M_2$ we have the following two commuting left $\GL_2$-actions $\rho_1, \rho_2: \GL_2 \to \GL(M_2)$ given by
\[
\rho_1(A)(B) \colon = AB \qquad \rho_2(A)(B) \colon = B \Abar
\]
for $A\in \GL_2$ and $B\in M_2$. Here for $A = \mat{a}{b}{c}{d}\in M_2$, $\Abar$ denotes the matrix $\mat{d}{-b}{-c}{a}$. Let $\cE(M_2)$ be the convergent filtered $F$-isocrystal attached to $(M_2, \rho_1, \rho_2)$. 

\begin{lemma}
\label{lemma:h1dr}
There exists an (up to scaling) canonical isomorphism  of filtered isocrystals on $\cH_p$ 
\[
\bH^1_{DR}(\cA/X_M) \cong \cE(M_2).
\]
\end{lemma}

{\em Proof.} This is essentially proved in section 5 of \cite{faltings2}. We will explain how the result shown there can be reformulated as above.

Let $D = \cB_p$ be the unique quaternion algebra over $\Q_p$ with maximal order $\cO_D$. Explicitely it can be written as $D = \Q_{p^2}[\Pi]$ with $\Pi^2 = p, \Pi x = \s(x) \Pi$. We recall the description of $\wh{\cH}$ as the moduli space of special formal $\cO_D$-modules. Let $R$ be a $\Zur$-algebra. A formal group $G$ over $R$ with $\cO_D$-multiplication is called special if the tangent space of $G$ is a free $\Z_{p^2}\otimes_{\Z_p} R$-module of rank 1 (where $\Z_{p^2}$ acts on $G$ via $\Z_{p^2}\hookrightarrow D$). Fix a special formal group $G_0$ over $\Fpbar$. Then $\wh{\cH}$ represents the functor which associates to every $\Zur$-algebra $R$ on which $p$ is nilpotent the set of isomorphism classes of pairs $(G, \la)$ where $G$ is a special formal group over $R$ and $\la: \Gbar \colon = G\otimes_R R/pR \to G_0\otimes_\Fpbar R/pR$ is a $\Q_p$-isogeny of height 0 of special formal $\cO_D$ -modules, i.e.\ $\la = p^{-h} \la_0$ for some positive integer $h$ and isogeny $\la_0$ of height $h$. Let $(\cG, \la)$ be the universal such pair over $\wh{\cH}$. The algebra $\End_{\cO_D}(G_0)_p$ of $\Q_p$-homomorphisms of $G_0$ is isomorphic to $M_2(\Q_p)$. We fix such an isomorphism $\End_{\cO_D}(G_0)_p \cong M_2(\Q_p)$. It induces compatible $\GL_2(\Q_p)$-actions on $\cG$ and $\wh{\cH}$.

Let $H^1_{\DR}(\cG/\wh{\cH})$ be the dual of the Lie algebra of the universal vectorial extension of $\cG$. It is a filtered convergent $F$-isocrystal on $\wh{\cH}$. Let $\DD(G_0)$ be the Dieudonn{\'e} module of $G_0$ and define $H^1_{\cris}(G_0)$ to be the dual of $\DD(G_0)$. Then as a convergent $F$-isocrystal $H^1_{\DR}(\cG/\wh{\cH})$ is constant isomorphic to $H^1_{\cris}(G_0) \otimes \cO_{\cH_p}$ (the isomorphism is induced by the isogeny $\la: \overline{\cG} \to G_0$). The $\Qur$-vector space $H^1_{\cris}(G_0)$ is endowed with the following structure:

\noi $\bu$ A $\s$-linear Frobenius $\Phi: H^1_{\cris}(G_0) \to H^1_{\cris}(G_0)$ such that $(H^1_{\cris}(G_0), \Phi)$ is isotypical of slope $\frac{1}{2}$.

\noi $\bu$ A $D$-module structure given by an embedding $j: D \to \End_{\Qur}(H^1_{\cris}(G_0))$ which commutes with $\Phi$ and such that $H^1_{\cris}(G_0)$ is a free $D\otimes_{\Q_p} \Qur$-module of rank 1.

\noi $\bu$ A $\Q_p$-algebra embedding $\io: M_2(\Q_p) \lra \End_{\Qur}(H^1_{\cris}(G_0))$ (induced by the isomorphism $\End_{\cO_D}(G_0)_p \cong M_2(\Q_p)$) which commutes with the $D$-action and the Frobenius.

Set $\Phi' \colon = j(\Pi)^{-1} \Phi$. Then $(H^1_{\cris}(G_0), \Phi')$ is isotypical of slope 0 , hence $V \colon = H^1_{\cris}(G_0)^{\Phi'=id}$ is a four-dimensional $\Q_p$-vector space. Therefore $H^1_{\cris}(G_0) \break\cong V_{\Qur}$ and under this isomorphism $\Phi$ corresponds to $\phi_V \otimes \s$ where $\phi_V$ denotes the restriction of $j(\Pi)$ to $V$. Thus we obtain an isomorphism of convergent $F$-isocrystals
\begin{equation}
\label{eqn:faltings}
H^1_{\DR}(\cG/\wh{\cH}) \,\, \cong \,\, V\otimes_{\Q_p} \cO_{\wh{\cH}}
\end{equation}
where the Frobenius on the right hand side is given by $\phi_V \otimes \Phi_{\cO_{\cH_p}}$. By restriction, the embedding $\io$ induces an embedding $\io_1: M_2(\Q_p) \to \End_{\Q_p}(V)$ and it is shown in \cite{faltings2} that the filtration on $H^1_{\DR}(\cG/\wh{\cH})$ corresponds under (\ref{eqn:faltings}) to the filtration on $V\otimes_{\Q_p} \cO_{\wh{\cH}}$ induced by the representation $\eta_1 \colon = \io_1\mid_{\GL_2}: \GL_2 \to \GL(V)$ as in section \ref{sec:gltwo}.

Let $Z$ be the centralizer of $\io_1(M_2(\Q_p))$ in $\End_{\Q_p}(V)$. Note that $\p_V\in Z$ and $\p_V^2 = p$. Clearly $Z \cong M_2(\Q_p)$. Choose an embedding $\io_2: M_2(\Q_p) \to \End_{\Q_p}(V)$ with image $Z$ and such that $\io_2(\mat{0}{p}{1}{0}) = \p_V$. Let $\eta_2\colon = \io_2\mid_{\GL_2}: \GL_2 \to \GL(V)$. Then (\ref{eqn:faltings}) can be interpreted as an isomorphism $H^1_{\DR}(\cG/\wh{\cH}) \cong \cE((V, \eta_1, \eta_2))$. On the other hand $(V, \eta_1, \eta_2) \cong (M_2, \rho_1, \rho_2)$. Thus we finally obtain an isomorphism of filtered $F$-isocrystals
\begin{equation}
\label{eqn:faltings2}
H^1_{\DR}(\cG/\wh{\cH}) \,\, \cong \,\, \cE(M_2).
\end{equation}
and the Lemma follows by descending (\ref{eqn:faltings2}) to $X_{\G}$.
\enddemo

It is apparent from the construction of (\ref{eqn:faltings2}) that there is an isomorphism $D_{\Qur} = \cB_{\Qur} \cong M_2(\Qur)$ such that the $D_{\Qur}^*$-action on the left hand side of (\ref{eqn:faltings2}) is compatible with the $\rho_2$-action of $\GL_2(\Qur)$-on the right hand side.

Hence from Lemma \ref{lemma:h1dr} we obtain also an isomorphism of filtered isocrystals with $\cB_{\Qur}^*$-action 
\[
\bH^2_{DR}(\cA/X_M) \cong \Lambda^2 \bH^1_{DR}(\cA/X_M) \cong \cE(\Lambda^2 M_2)
\]
Note that $(M_2, \rho_1, \rho_2)$ is canonically isomorphic to $V_1\odot V_1$. Hence
\[
(\Lambda^2 M_2,\Lambda^2\rho_1,\Lambda^2\rho_2) = \Lambda^2 V_1 \odot \sym^2 V_1 \oplus \sym^2 V_1 \odot \Lambda^2 V_1 = \{1\}V_2 \oplus V_2\{1\}
\]
and therefore $\bigcap_{b\in \cB_{\Qur}^*} \Ker(b- N(b): \cE(\Lambda^2 M_2) \to \cE(\Lambda^2 M_2)) \cong \cV_2$.

Applying theorem \ref{theorem:comp} (a) we obtain
\[
D_{\st, \Qur}(H^1_{\et}(\overline{X}_M, \bL_2)) \cong H^1_{\DR}((X_M)_{\Qur}, \cV_2).
\]
Passing on both sides to $G$-invariants yields an isomorphism
\begin{equation}
\label{eqn:dieumodformn2}
D_{\st, \Qur}(H_p(\cM_2)) = H^1_{\DR}(X_{\G}, \cV_2)
\end{equation}
i.e.\ the assertion for $n=2$. 

Now let $n \geq 4$. There exists a canonical nondegenerated symmetric bilinearform in $\Rep_{\Q_p}(\GL_2)$ 
\begin{equation}
\label{eqn:s4pair2}
< \,\,\,,\,\,\,> : V_2 \otimes V_2 \lra {\det}^{\otimes^2}
\end{equation}
whose definition we briefly recall from (\cite{bdis}, 1.2; actually (\ref{eqn:s4pair2}) is $-2 \times$ the pairing considered there). Let $\cU = \{U \in M_2\mid \, \trace(U) = 0\}$ with right $\GL_2$-action 
\[
U\cdot A = \Abar\, U \, A
\]
for $U\in\cU$ and $A\in \GL_2$. For $U\in \cU$ we set
\begin{equation}
\label{eqn:defpu}
P_U(X) = \trace\left(U \mat{X}{-X^2}{1}{-X}\right) 
\end{equation}
The map $\cU \to \cP_2, U \mapsto P_U(X)$ is an isomorphism of right $\GL_2$-modules. There is a pairing on $\cP_2$ given by 
\[
< P_{U_1}(X), P_{U_2}(X) > \colon = - \trace(U_1 \Ubar_2)\qquad \forall U_1, U_2\in\cU.
\]
The isomorphism $\cP_2 \to V_2, R(X) \mapsto (P(X) \mapsto <R(X), P(X)>)$ yields the pairing (\ref{eqn:s4pair2}) by transport of structure. Note that $<A v_1, v_2> = < v_1, \Abar v_2>$ for all $v_1, v_2\in V_2$ and $A\in \GL_2$. 

By (\cite{bdis}, 1.2) the $\GL_2$-subrepresentation $V_{2m}$ of $\sym^m V_2$ is the kernel of the Laplace operator associated to (\ref{eqn:s4pair2})  
\begin{equation}
\label{eqn:lap2}
\De: \sym^m V_2 \lra \sym^{m-2} V_2\otimes {\det}^{\otimes^2}.
\end{equation}
By tensoring (\ref{eqn:s4pair2}) with $\odot \det^{\otimes^2}$ we obtain a pairing 
\begin{equation}
\label{eqn:s4pair3}
< \,\,\,,\,\,\,> : V_2\{1\} \otimes V_2\{1\} \lra \{2|2\}.
\end{equation}

We need the following elementary lemma.

\begin{lemma}
\label{lemma:pair} The pairing {\rm (\ref{eqn:s4pair3})} coincides with the pairing 
\[
V_2\{1\} \otimes V_2\{1\} \hookrightarrow \Lambda^2 M_2 \otimes \Lambda^2 M_2 \lra \Lambda^4 M_2 \cong \{2|2\} 
\]
(as a map in $\Rep_{\Q_p}(\GL_2\times \GL_2)$).
\end{lemma}

By applying theorem 1.5 (a) to $\cA^i \to X_M$ and using the Kuenneth formula  for $H^{2i}_{\DR}(\cA^i/X_M)$ one can easily deduce from Lemma \ref{lemma:h1dr} that
\[
D_{\st, \Qur}(H^1(\overline{X}_M, \sym^i\bL_2)) \cong H^1_{\DR}((X_M)_{\Qur}, \sym^i\cV_2),
\]
for all $i\geq 1$. Lemma \ref{lemma:pair} implies that the diagram 
\[
\begin{CD}
D_{\st, \Qur}(H^1(\overline{X}_M, \sym^m\bL_2)) @> \cong >> H^1_{\DR}((X_M)_{\Qur}, \sym^m\cV_2)\\
@VVV@VVV\\
D_{\st, \Qur}(H^1(\overline{X}_M, \sym^{m-2}\bL_2)) @> \cong >> H^1_{\DR}((X_M)_{\Qur}, \sym^{m-2}\cV_2),
\end{CD}
\]
commutes where the vertical maps are induced by (\ref{eqn:lap1}) and (\ref{eqn:lap2}) respectively. Hence the kernels of the vertical arrows are isomorphic as well. Passing again to $G$-invariants yields Theorem \ref{theorem:dieumodform} in the case $n\geq 4$.\enddemo

\begin{remark} 
\label{remark:dieumodformpair}
{\rm The pairing (\ref{eqn:s4pair3}) induces a symmetric bilinearform 
\begin{equation}
\label{eqn:s4pair4}
< \,\,\,,\,\,\,> : V_n \otimes V_n \lra \det^{\otimes^n}
\end{equation}
(compare \cite{bdis}, 1.2). It follows from Lemma \ref{lemma:pair} that the (nondegenerate) pairing
\begin{equation}
\label{eqn:s4pair5}
< \,\,\,,\,\,\,> : H^1_{\DR}(X_{\G}, \cV_n)\otimes H^1_{\DR}(X_{\G}, \cV_n) \lra H^2_{\DR}(X_{\G}, \cV_n\otimes \cV_n) \lra
\end{equation}
\[
\stackrel{(*)}{\lra} H^2_{\DR}(X_{\G}) \cong \Qur
\]
(where $(*)$ is induced by (\ref{eqn:s4pair4})) is given by the cup-product
\[
H^{n+1}_{\DR}(\cM_n) \otimes H^{n+1}_{\DR}(\cM_n) \seq H^{n+1}_{\DR}(\cA^m) \otimes H^{n+1}_{\DR}(\cA^m)\lra H^{2n+2}_{\DR}(\cA^m) \cong \Qur,
\]
hence under the isomorphism (\ref{eqn:dieumodformn2}) coincides with the map
\begin{equation}
\label{eqn:s4pair6}
D_{\st, \Qur}(H_p(\cM_n)) \otimes D_{\st, \Qur}(H_p(\cM_n)) \lra \Qur[2m+1]
\end{equation}
induced by the cup-product $\bL_n \otimes \bL_n \to \Q_p[-n]$. More precisely the isomorphism (\ref{eqn:dieumodformn2}) can and will be chosen in such a way that it is compatible with to (\ref{eqn:s4pair5}) and (\ref{eqn:s4pair6}). This fixes it up to sign (and not only up to scaling).

Hence we can view (\ref{eqn:s4pair5}) as a map
\begin{equation}
\label{eqn:s4pair7}
< \,\,\,,\,\,\,> : H^1_{\DR}(X_{\G}, \cV_n)\otimes H^1_{\DR}(X_{\G}, \cV_n) \lra \Qur[2m+1]
\end{equation}
in $MF_{\Qur}(\p,N)$.}
\end{remark}

We finish this section by describing explicitely the sequence of Dieudonn\'e modules of a certain Gysin sequence associated to a finite number of rational points on $X$. Let $x_1, \ldots, x_r\in X(\Qur)$ and fix points $z_1, \ldots, z_r \in \cH_p(\Qur)$ such that $\G z_i = x_i$ (under the isomorphism $X(\Qur) \cong \G \backslash \cH_p(\Qur)$. For simplicity we assume that the stabilizers of $z_1, \ldots, z_r$ in $\G$ are all $\{\pm 1\}$. Let $X_0 \cong \G_M \backslash \cH_p$ be one of the components of $(X_M)_{\Q_p}$ and let $x'_i \colon = \G_M z_i\in X'$ and $U'$ their complement in $(X_M)_{\Q_p}$. We have a short exact (Gysin) sequence of $G_{\Qur}$-modules
\[
0 \lra H^1(\overline{X}_M, \bL_n) \lra H^1(U\otimes_{\Qur} \overline{\Qur}, \bL_n) \lra \bigoplus_{i=1}^r\, (\bL_n)_{\bar{x}'_i}(-1) \lra 0.
\]
Let 
\begin{equation}
\label{eqn:s4gysin}
0 \lra H_p(\cM_n) \lra E \lra \bigoplus_{i=1}^r\, (\bL_n)_{\bar{x}'_i}(-1) \lra 0
\end{equation}
be its push-out under the canonical projection $H^1(\overline{X}_M, \bL_n) \to H_p(\cM_n)$. Similarly as above we can deduce from Theorem \ref{theorem:comp} (b) the following

\begin{theorem}
\label{theorem:dieumodformgys} 
After applying $D_{\st, \Qur}$ the sequence {\rm (\ref{eqn:s4gysin})} becomes isomorphic to the sequence of filtered $(\p,N)$-modules
\[
\begin{CD}
0 @>>> H^1_{\DR}(X_{\G}, \cV_n) @>>> H^1_{\DR}(U, \cV_n) @> \oplus_{i}\,\Res_{z_i} >> \bigoplus \, (\cV_n)_{z_i}[1] @>>> 0
\end{CD}
\]
where $U \colon = X_{\Qur}- \{x_1, \ldots, x_r\}$.
\end{theorem}

\begin{remark} 
\label{remark:dieumodformpair2}
{\rm In particular for $x\in X_M(\Qur)$ and $z\in \cH_p$ lying above $x$ there is a canonical isomorphism
\[
D_{st}((\bL_n)_{\bar{x}}) \cong (\cV_n)_z.
\]
It is easy to see that it respects the pairings on both sides induced by the cup product $\bL_n \otimes \bL_n \to \Q_p[-n]$ and (\ref{eqn:s4pair4}) respectively.}
\end{remark}

\section{Comparison of $\cL$-invariants}
\label{sec:linv} 

As in the last section let $f_{\infty}$ be a newform of weight $k\geq 4$ and level $N$ corresponding to a newform $f$ on $X$. The exceptional zero conjecture  of \cite{mazur_tate_teitelbaum} relates the derivative of a $p$-adic $L$-function to the value of the complex $L$-function of $f_{\infty}$ at the central critical point $\frac{k}{2}$. In this conjecture a certain local factor $\cL(f_{\infty})$ appears the so called $\cL$-invariant of $f_{\infty}$. Three possible definitions for it have been given. The Fontaine-Mazur $\cL$-invariant $\cL_{FM}(f_{\infty})$ is defined in terms of the semistable Dieudonn\'e module $D_{\st}(V_p(f_{\infty}))$. The $\cL$-invariants of Teitelbaum and Coleman, $\cL_{T}(f_{\infty})$ and $\cL_{C}(f_{\infty})$ are defined in terms of $p$-adic integration on the Shimura curve $X$ and the modular curve $X_0(N)$ respectively. The exceptional zero conjecture for $f_{\infty}$ has been proved by Stevens (using $\cL_{C}(f_{\infty})$) and by Kato, Kurihara and Tsuji (using $\cL_{FM}(f_{\infty})$). In \cite{bdis} a version of the exceptional zero conjecture is proved involving the anticyclotomic $p$-adic $L$-function of $f_{\infty}$ and $\cL_{T}(f_{\infty})$.

In this section we show that $\cL_{FM}(f_{\infty}) = \cL_{T}(f_{\infty})$. Our proof is based on the explicit description given in the previous section of the Dieudonn\'e module of $V_p(f_{\infty}) \cong V_p(f)$ (or more precisely of $H_p(\cM_n)$) considered as a representation of the inertia group $I_p \cong G_{\Qur}$ which in turn is based on the comparison theorems of Faltings and Coleman-Iovita. 

We begin by giving an explicite description of the Hodge filtration of $H^1(X_{\G}, \cV_n) = H^1_{\DR}(X_{\G}, \cV_n)$. 

\begin{proposition}
\label{proposition:filhodge} We have,
\[
F^i H^1(X_{\G}, \cV_n) =  
\left\{ \begin{array}{ll}
                  H^1(X_{\G}, \cV_n) & \mbox{if $i\leq 0$},\\
                  M_k(\G) & \mbox{if $1 \leq i\leq k-1$},\\
                  0 & \mbox{if $i\geq k$}
\end{array}\right.
\]
\end{proposition}

{\em Proof.} For $i\in \{1, \ldots, n\}$ we let $\partial^i\in V_n\otimes \cO(\cH_p)$ be given by $\partial^i(z)(P(X)) \break \colon = (\frac{d^i}{dX^i} P(X))|_{X=z}$ for $z\in \cH_p, P(X) \in \cP_n$. It is easy to verify (and will be left to the reader) that the filtration $\cF^{\bu} \cV_n$ on $\cV_n$ is given by 
\[
\cF^j \cV_n = \left\{ \begin{array}{ll}
                  \cV_n & \mbox{if $j\leq 0$},\\
                  \sum_{i=0}^{n-j}\, \cO_{\cH_p}\, \partial^i & \mbox{if $0 \le                  j\leq n$},\\
                  0 & \mbox{if $i\geq n+1$}
\end{array}\right.
\]
Consequently $F^i H^1(X_{\G}, \cV_n) = H^1(X_{\G}, \cV_n)$ if $i\leq 0$, $F^i H^1(X_{\G}, \cV_n) = 0$ if $i\geq n+2$ and the image of the embedding
\[
M_k(\G) \lra H^1(X_{\G}, \cV_n), f(z) \mapsto \omega_f \colon = f(z) \partial^0 \otimes dz
\]
lies in $F^{n+1} H^1(X_{\G}, \cV_n)$. Hence,
\[
\dim F^{n+1} H^1(X_{\G}, \cV_n) \geq \dim M_k(\G) = \frac{1}{2} \dim H^1(X_{\G}, \cV_n)
\]
(see \cite{shalit3}). Since $F^1$ and $F^{n+1}$ are orthogonal with respect to (\ref{eqn:s4pair7}) we obtain
\[
\dim F^{n+1} \leq \dim F^1 \leq \dim H^1(X_{\G}, \cV_n)/F^{n+1} \leq \frac{1}{2} \dim H^1(X_{\G}, \cV_n)
\]
and therefore $F^1 = F^2 = \ldots = F^{n+1} = M_k(\G)$. \enddemo

By Theorem \ref{theorem:dieumodform}, $V_p(f)$ -- considered as a $G_{\Q_p}$-representation -- is semisimple. The associated Dieudonn{\'e} module $D_{\st, \Q_p}(V_p(f))$ is a two-dimensional $F_{\Q_p}$-module where $F = F_f$. The Fontaine-Mazur $\cL$-invariant of $f_{\infty}$ is defined as the $\cL$-invariant of $D_{\st, \Q_p}(V_p(f_{\infty})) \cong D_{\st, \Q_p}(V_p(f))$. The following result has been conjectured by Mazur (\cite{mazur}, section 12).

\begin{lemma}
\label{lemma:cmd} $D_{\st, \Q_p}(V_p(f))$ is a monodromy $F_{\Q_p}$-module.
\end{lemma}

{\em Proof.} By remark \ref{remark:Lbasechange} (a) above it is enough to consider $D_{\st, \Qur}(V_p(f))$ instead. Clearly property (i) of \ref{definition:monodmod} is satisfied and that (ii) holds is a consequence of the exactness of (\ref{equation:NNex}). Proposition \ref{proposition:filhodge} implies that $F^{m+1} D_{\st, \Qur}(V_p(f))\break = e_f(F^{m+1} H^1(X_{\G}, \cV_n)) \cong e_f M_k(\G)$, hence $F^{m+1} D_{\st, \Qur}(V_p(f))$ is free a $F_{\Qur}$-module of rank 1. Finally the fact that the restriction of $N:H^1(X_{\G}, \cV_n)\to H^1(X_{\G}, \cV_n)$ to $M_k(\G)$ is injective (see \cite{shalit1}, 3.9) implies that 
\[
F^{m+1} D_{st,\Qur}(V_p(f)) \cap \Ker(N) =0.
\]
\enddemo

\begin{definition}
\label{definition:fmcLinv}
The Fontaine-Mazur $\cL$-invariant $\cL_{FM}(f)$ of $f$ (or $f_{\infty}$) is defined as the $\cL$-invariant of the monodromy module $D_{st,\Q_p}(V_p(f))$.
\end{definition}

We also recall the definition of the Teitelbaum $\cL$-invariant $\cL_T(f)$. Let $P_0: M_k(\G) \to H^1(\G, V_n)$ be the composition of the inclusion
\[
M_k(\G)\cong F^{m+1} H^1(X_{\G}, \cV_n) \hookrightarrow  H^1(X_{\G}, \cV_n)
\]
with the map $P$, i.e.\ $P(f)$ is represented by the cocycle $\ga \mapsto \ga(F_f) - F_f$ where $F_f$ is a primitive of $\omega_f$. Both maps $P_0$ and
\begin{equation}
\label{eqn:ischneider}
M_k(\G) \hookrightarrow H^1(X_{\G}, \cV_n) \stackrel{\ep\circ I}{\lra} H^1(\G, V_n)
\end{equation} 
are homomorphisms of free $\TT_{\Qur}$-modules of rank one and (\ref{eqn:ischneider}) is an isomorphism. Hence there is an element $\cL_T\in \TT_{\Qur}$ with $P_0(g) = \cL_T \ep(I(\omega_g))$ for every $g\in M_k(\G)$. Then $\cL_T(f) \colon = e_f \cL_T\in e_f\TT_{\Qur} = F_{\Qur}$.

\begin{theorem} 
\label{theorem:lfm=lt}
Let $f\in M_k(X, \Q)$ be a newform. Then, $\cL_{FM}(f) = \cL_T(f)$.
\end{theorem}

{\em Proof.} To simplify the notation we work with $D_{\st}(H_p(\cM_n)) \cong H^1(X_{\G}, \cV_n)$ rather than $D_{\st, \Qur}(V_p(f))$ (though strictly speaking $D_{\st, \Qur}(V_p(f))$ is not a monodromy $\TT_{\Q_p}$-module in the sense of definition \ref{definition:monodmod} since $\TT_{\Q_p}$ is in general not semisimple). The slope decomposition of $H^1(X_{\G}, \cV_n)$ is of the form 
\begin{equation}
\label{eqn:s5slope1}
H^1(X_{\G}, \cV_n) = H^1(X_{\G}, \cV_n)_m \oplus H^1(X_{\G}, \cV_n)_{m+1}.
\end{equation}
where $H^1(X_{\G}, \cV_n)_m = \io(H^1(\G,V_n))= \Ker(N)$ and $H^1(X_{\G}, \cV_n)_{m+1}$ is the kernel of $P$. If we decompose an element $x\in H^1(X_{\G}, \cV_n)$ as $x= x_m + x_{m+1}$ according to (\ref{eqn:s5slope1}) then $\cL_T$ is thus characterised by the property 
\[
\cL_T N(x) = - x_m\qquad \forall \, x\in F^{m+1} H^1(X_{\G}, \cV_n),
\]
whereas $\cL_{FM} = \cL_{FM}(H^1(X_{\G}, \cV_n))\in \TT_{\Qur}$ is characterised by
\[
x - \cL_{FM} N(x)\in F^{m+1} H^1(X_{\G}, \cV_n)\qquad\forall \, x\in H^1(X_{\G}, \cV_n)_{m+1}.
\]
By (\cite{shalit3}, Theorem 1.6) any element $x$ of $H^1(X_{\G}, \cV_n)$ can be uniquely written as $x = x' + x''$ with $x'\in F^{m+1}$ and $x''\in H^1(X_{\G}, \cV_n)_m$. Therefore if $x\in H^1(X_{\G}, \cV_n)_{m+1}$ then $x'_m = - x''$ and thus
\[
x - \cL_T N(x) =  x' + x'' - \cL_T N(x') = x' + x'' + x'_m = x' \in F^{m+1}.
\]
Since this holds for all $x\in H^1(X_{\G}, \cV_n)_{m+1}$ we conclude $\cL_T = \cL_{FM}$.
\enddemo

\section{The $p$-adic Abel-Jacobi map}
\label{sec:pabeljac}

In this section we will define the $p$-adic Abel-Jacobi maps for the motives $\cM_n$ (see appendix \ref{subsection:a1}) where $n=2m$ is a positive even integer. It is a map from the Chow group $CH^{m+1}(\cM_n)$ to the dual of the space of weight $k$-modular modular forms on $X$.

We begin by reviewing briefly the definition of the cohomological $\ell$-adic Abel-Jacobi map (see e.g.\ \cite{jannsen}). Let $K$ be a field of characteristic 0 and $\ell$ a prime number. For a smooth projective variety $X$ over $K$ we denote by $CH^i(X)= CH^i(X)_{\Q}$ the Chow group of codimension $i$-cycles (with rational coefficients) and by $CH^i(X)_0$ be the subgroup of cycle classes homologous to zero, i.e. the kernel of the cycle class map
\begin{equation}
\label{eqn:cl}
cl=cl^{X,i}:  CH^i(X) \lra H^{2i}(\Xbar, \Ql(i))^{G_K}
\end{equation}
where $\Xbar \colon = X\otimes_K \Kbar$. The map (\ref{eqn:cl}) factors through $H^{2i}_{\cont}(X, \Ql(i))$. The $\ell$-adic Abel-Jacobi map 
\begin{equation}
\label{eqn:aj1}
cl_0 = cl^{X,i}_0:  CH^i(X)_0 \lra H^1_{\cont}(K, H^{2i-1}(\Xbar, \Ql(i))) = \Ext^1_{G_K}(\Ql, H^{2i-1}(\Xbar, \Ql(i)))
\end{equation}
is defined as follows: Let $z\in CH^i(X)_0$ and $Z$ be a cycle representing $z$. Then $cl_0(z)$ is the extension class
\begin{equation}
\label{eqn:aj2}
0 \lra H^{2i-1}(\Xbar, \Ql(i)) \lra E \lra \Ql \lra 0
\end{equation}
given by the pull-back of
\[
0 \lra H^{2i-1}(\Xbar, \Ql(i)) \lra H^{2i-1}(\Xbar- |\Zbar|, \Ql(i)) \lra 
\]
\[
\lra\Ker(H^{2i}_{|\Zbar|}(\Xbar, \Ql(i))\to H^{2i}(\Xbar, \Ql(i))) \lra 0
\]
via $\Ql\cong \Ql\cdot cl^{\Xbar}_{\Zbar}(\Zbar)\hookrightarrow H^{2i}_{|\Zbar|}(\Xbar, \Ql(i))$.\\

Let $K$ be a finite extension of $\Q_p$. For any $p$-adic Galois representation $V$ of $G_K$, Bloch and Kato and Nekovar have introduced subspaces
$H^1_f(K, V)\seq H^1_{\st}(K, V) \seq H^1_g(K, V)$ of $H^1(K, V)$ given by
\[
\begin{array}{c}
H^1_f(K, V) = \Ker(H^1(K, V) \to H^1(K, V\otimes_{\Q_p} B_{\cris})\\
H^1_{\st}(K, V) = \Ker(H^1(K, V) \to H^1(K, V\otimes_{\Q_p} B_{\st})\\
H^1_g(K, V) = \Ker(H^1(K, V) \to H^1(K, V\otimes_{\Q_p} B_{\DR})
\end{array}
\]
(for the first and the last see \cite{bloch_kato} and for the second see \cite{nekovar2}). It is known \cite{nekovar2} that if $V$ is a semistable (resp.\ cristalline, resp.\ de Rham) representation of $G_K$ then $H^1_{\st}(K, V)$ (resp.\ $H^1_f(K, V)$, resp.\ $H^1_g(K, V)$) can be identified with the group of extension classes of $V$ by $\Q_p$ in the category of semistable (resp.\ cristalline resp.\ de Rham) representations of $G_K$. 

\begin{lemma}
\label{lemma:nekovar}
Let $X$ be smooth projective variety over $K$. Then the image of (\ref{eqn:aj1}) for $\ell = p$ is contained in $H^1_{\st}(K, H^{2i-1}(\Xbar, \Q_p(i)))$.
\end{lemma}

Nekovar has anounced a proof of this fact in (see \cite{nekovar4}). His argument is based on a spectral sequence relating $\Ext$-groups in the category of semistable $p$-adic representations to log-syntomic cohomology. If $X$ has good reduction a proof along these lines appears already in \cite{nekovar5}. Nekovar has indicated also a different approach for proving this lemma (\cite{nekovar4}, 3.11): firstly by using de Jong's theorem on alterations he shows that it is enough to consider the case where $X$ has semistable reduction (\cite{nekovar4}, page 6 above). Then $H^{2i-1}(\Xbar, \Q_p(i)))$ is a semistable representation (see \cite{tsuji}) and therefore $H^1_{\st}(K, H^{2i-1}(\Xbar, \Q_p(i))) = H^1_g(K, H^{2i-1}(\Xbar, \Q_p(i)))$ by (\cite{hyodo} or \cite{nekovar2}, proposition 1.24). Consequently it is enough to show that in the extension (\ref{eqn:aj2}) corresponding to an element $z\in CH^i(X)_0$, the middle term is a de Rham representation or -- if $Z$ be a cycle representing $z$ -- that $H^{2i-1}(\Xbar- |\Zbar|, \Q_p(i))$ is de Rham. However this follows from \cite{kisin} where it is shown that $H^{2i-1}(\Xbar-|\Zbar|, \Q_p(i))$ is potentially semistable hence de Rham.\\

Now assume that $K$ is a number field and let $v$ be a place of $K$ above $p$ which for simplicity we assume to be unramified over $p$. Let $X$ be a smooth projective variety over $K$ and assume that $H^j(\Xbar, \Q_p)$ is semistable as a representation of the local Galois group $G_{K_v}$ for all $j$. Then the functor $D_{\st, K_v}$ yields isomorphisms
\begin{equation}
\label{eqn:isoh0}
 H^{2i}(\Xbar, \Q_p(i))^{G_{K_v}} \cong \G(D_{\st, K_v}(H^{2i}(\Xbar, \Q_p(i))),
\end{equation}
\begin{equation}
\label{eqn:isoh1}
H^1_{\st}(K_v, H^{2i-1}(\Xbar, \Q_p(i))) \cong \Ext^1_{\sRep{K_v}}(\Q_p,H^{2i-1}(\Xbar, \Q_p(i))) \cong 
\end{equation}
\[
\cong \Ext^1_{MF^{ad}_{K_v}(\p, N)}(K_v[i], D_{\st}(H^{2i-1}(\Xbar, \Q_p))),
\]
where for a filtered $(\p,N)$-module $D$, $\G(D)$ is given by
\[
\G(D) = \Hom_{MF^{ad}_{K_v}}(K_v, D) = F^0 D \cap D^{\p = id, N =0}.
\]
By composing (\ref{eqn:cl}) and (\ref{eqn:aj1}) for $\ell =p$ with the restriction maps for $K_v/K$ and with (\ref{eqn:isoh0}) and (\ref{eqn:isoh1}) respectively we obtain maps
\begin{equation}
\label{eqn:cl2}
cl=cl^{X,i}:  CH^i(X) \lra \G(D_{\st, K_v}(H^{2i}(\Xbar, \Q_p(i)))),
\end{equation}
\begin{equation}
\label{eqn:aj3}
cl_0 = cl^{X,i}_0:  CH^i(X)_0 \lra \Ext^1_{MF^{ad}_{K_v}(\p, N)}(K_v[i], D_{\st, K_v}(H^{2i-1}(\Xbar, \Q_p))).
\end{equation}

One can easily extend the definition of (\ref{eqn:cl}), (\ref{eqn:aj1}), (\ref{eqn:cl2}) and (\ref{eqn:aj3}) to Chow motives. In the following we use the notation of the appendix \ref{subsection:a1} (for definitions and general facts about Chow motives we refer to \cite{scholl2}). Note that $\Corr(X, X)$ acts on the source and the target of (\ref{eqn:cl}) and (\ref{eqn:aj1}) and that both maps are homomorphisms of $\Corr(X, X)$-modules. Thus if $\cM =(X,p)$ is a motive over $K$ such that its $p$-adic realizations are semistable as $G_{K_v}$-representations we obtain maps
\begin{equation}
\label{eqn:clm}
cl = cl^{\cM, i}:  CH^i(\cM) \lra \G(D_{\st, K_v}(H^{2i}_p(\cM)(i))),
\end{equation}
\begin{equation}
\label{eqn:aj4}
CH^i(\cM)_0 \lra \Ext^1_{MF^{ad}_{K_v}(\p, N)}(K_v[i], D_{\st, K_v}(H^{2i-1}_p(\cM))).
\end{equation}

This applies in particular to the sequence of motives $\cM_n$. By Lemma \ref{lemma:a1hpm} of the appendix we see that $CH^{m+1}((\cM_n)_K)_0 = CH^{m+1}((\cM_n)_K)$. Lemma \ref{lemma:bkexp} together with (\cite{nekovar2}, 1.27) (and appendix \ref{subsection:a1}) shows that the target of (\ref{eqn:aj4}) for $\cM = (\cM_n)_K$, $i=m+1$ can be identified with
\begin{equation}
\label{eqn:drm}
\Ext^1_{MF^{ad}_{K_v}(\p, N)}(K_v[m+1], D_{\st, K_v}(H_p(\cM_n))) \cong H^{n+1}_{\DR}((\cM_n)_{K_v})/F^{m+1}
\end{equation}
\[
\cong F^{m+1} H^{n+1}_{\DR}((\cM_n)_{K_v})^{\vee} \cong M_k(X, K_v)^{\vee}.
\]

\begin{definition}
\label{definition:paj}
Let
\[
\rho_K = \rho_{K,v}: CH^{m+1}((\cM_n)_K)\lra M_k(X, K_v)^{\vee}
\]
be the composite of {\rm (\ref{eqn:aj4})} (for $\cM = \cM_n$, $i=m+1$) with {\rm(\ref{eqn:drm})}. It will be called the $v$-adic Abel-Jacobi map for $(\cM_n)_K$.
\end{definition}

Let $K \hookrightarrow \Qur$ be an inclusion which induces the place $v$ of $K$. By abuse of notation we will denote the composition of $\rho_{K,v}$ with the inclusion $M_k(X, K_v)^{\vee} \hookrightarrow M_k(X, \Qur)^{\vee} \cong M_k(\G)^{\vee}$ also by $\rho_K$
\begin{equation}
\label{eqn:s7paj}
\rho_K = \rho_{K,v}: CH^{m+1}((\cM_n)_K)\lra M_k(\G)^{\vee}.
\end{equation}

Let $\wt{P}$ be a closed point of $(X_M)_K$ and let $\cA^m_{\wt{P}}$ be the fiber of $(\cA^m)_K \to (X_M)_K$ over $\wt{P}$. We assume that the residue field $L \colon = K(\wt{P})$ of $\wt{P}$ is Galois over $\Q$ and that $p$ is unramified in $L$. We let $P$ be the image of $\wt{P}$ under $p: (X_M)_K \to X_K$ with residue field $H$ (later $K$ will be an imaginary quadratic field and $P$ a Heegner point so that $H$ is the Hilbert class field of $K$ and $L = H(\mu_M)$). Our aim is to compute the image of a cycle class under (\ref{eqn:s7paj}) which has support on $\cA^m_{\wt{P}}$. More precisely we fix an embedding $L \to \Qur$ which induces the place $v$ on $K$ (the induced place on $H$ and $L$ will be also denoted by $v$). It determines a point on $(X_M)_{\Qur}$ above $\wt{P}$ which according to (\ref{eqn:puniv2}) can be written as $\G_M z, z\in \cH_p$ (lying in one of the components $\G_M\backslash \cH_p$). We also assume that the stabilizer of $z$ (in $\G$) is $\{\pm 1\}$. 

There is an exact sequence of $G_L$-modules 
\begin{equation}
\label{eqn:aj5}
0 \lra H^1(\overline{X}_M, \bL_n)(m+1) \lra H^1(\overline{X}_M - (\wt{P}\otimes_L \Qbar), \bL_n)(m+1) \lra
\end{equation}
\[
\lra (\bL_n)_{\wt{P}\otimes_L \Qbar}(m) = H^{2m}_p((\cA^m_{\wt{P}}, \ep_n)/L)(m)\lra 0,
\]
where $H^{2m}_p((\cA^m_{\wt{P}}, \ep_n)/K)$ is the $p$-adic realisation of the motive $(\cA^m_{\wt{P}}, \ep_n)/L$. Let 
\begin{equation}
\label{eqn:aj6}
0 \lra H_p(\cM_n)(m+1) \lra E \lra (\bL_n)_{\wt{P}\otimes_L \Qbar}(m)\lra 0
\end{equation}
be the push-out of (\ref{eqn:aj5}) under $H^1(\overline{X}_M, \bL_n) \to H_p(\cM_n)$. By restricting the Galois action to $G_{\Qur}\seq G_L$ and applying $D_{\st, \Qur}$ to (\ref{eqn:aj6}) we obtain a short exact sequence in $MF^{ad}_{\Qur}(\p, N)$ which according to Theorem \ref{theorem:dieumodformgys} is isomorphic to the sequence
\begin{equation}
\label{eqn:gysin4}
0 \lra H^1(X_{\Qur}, \cV_n)[-(m+1)] \lra H^1(U, \cV_n)[-(m+1)] \stackrel{\Res_z}{\lra} (\cV_n)_z[-m] \lra 0
\end{equation}
where $U \colon = X_{\G}-\{\G z\}$.

\begin{proposition} 
\label{prop:ajfiber}
The diagram 
\begin{equation}
\label{eqn:cdaj1}
\begin{CD}
CH^m((\cA^m_{\wt{P}}, \ep_n)) @>>> \G((\cV_n)_z[-m])\\
@VVV@VVV\\
CH^{m+1}((\cM_n)_H) @>\rho_H >> M_k(\G)^{\vee}
\end{CD}
\end{equation}
commutes. Here the upper horizontal map is the composition of {\rm (\ref{eqn:clm})} for $\cM = (\cA^m_{\wt{P}}, \ep_n)/L$ with the restriction map to $\Qur$.
The right vertical map is induced by $\cA^m_{\wt{P}} \hookrightarrow (\cA^m)_H$ and the left is a connecting homomorphism in the long exact $\Ext$-groups sequence corresponding to {\rm (\ref{eqn:gysin4})}.
\end{proposition}

{\em Proof.} This can be easily deduced from the commutativity of the diagram 
\begin{equation}
\label{eqn:cdaj2}
\begin{CD}
\ep_n CH^m((\cA^m_{\wt{P}}, \ep_n)) @>>> \ep_n H^{2m}(\cA^m_{\wt{P}}\otimes_L\Qbar, \Q_p(m))^{G_L}\\
@VVV@VV \pa V\\
\ep_n CH^{m+1}((\cA^m)_H) @>>> H^1_{\cont}(L, \ep_n H^{2m-1}((\cA^m)_\Qbar, \Q_p(m)))
\end{CD}
\end{equation}
where the lower arrow is given by the composition 
\[
\ep_n CH^{m+1}((\cA^m)_H)\lra H^1_{\cont}(H, \ep_n H^{2m-1}((\cA^m)_\Qbar, \Q_p(m)))\lra 
\]
\[
\stackrel{res}{\lra} H^1_{\cont}(L, \ep_n H^{2m-1}((\cA^m)_\Qbar, \Q_p(m))).
\]
\enddemo

\section{Heegner cycles}
\label{sec:heegner}

In this section we define Heegner cycles. They are not really algebraic cycles but elements in $CH^{m+1}(\cM_n)$. We then compute their images under the $p$-adic Abel-Jacobi map.

To begin with we recall the definition of CM-points and Heegner points on $X = X_{N^+, pN^{-}}$. Let $F$ be a field of characteristic 0. Recall that any $F$-valued point of $X$ can be represented by an abelian surface with QM, i.e.\ by a triple $(A, \iota, C)$ where $A$ is an abelian surface, $\iota: \cR^{\max} \to \End_F(A)$ a ring monomorphism and $C$ a cyclic subgroup of order $N^+$ in $A$ which is stable under the action of the Eichler order $\cR$. The abelian surface $A$ is said to have complex multiplication if $\End_{\cR^{\max}}(A) \neq \Z$. In this case it is an order $\cO$ in an imaginary quadratic number field and $A$ is said to have complex multiplication by $\cO$. An $F$-valued point on $X$ is called CM-point if it can be represented by a triple $(A, \iota, C)$ such that $A$ has complex multiplication. It is called Heegner point if $(A, \iota, C)$ can be chosen so that $A$ has complex multiplication by $\cO$ and $C$ is $\cO$-stable.

Before we proceed with the definition of Heegner cycles we need to determine the structure of the Neron-Severi group (tensored by $\Q$) of an abelian surface with QM. Let $A$ be an abelian surface over $F$, $\io: \cR^{\max} \to \End_F(A)$ a monomorphism and set $\NS(A)_{\Q} = \NS(A)\otimes_{\Z} \Q$. We have a natural right $\cB^*$-action on $\NS(A)_{\Q}$ given by $\cL\cdot b = \io(b)^*(\cL)$. By $\ad^0(\cB)$ we denote the representation of $\cB^*$ consisting of elements of trace $=0$ on which $\cB^*$ acts from the right by $u\cdot b = \bar{b} u b$. 

\begin{lemma}
\label{lemma:l5} As a right $\cB^*$-module $\NS(A)\otimes \Q$ decomposes as follows into irreducible $\cB^*$-representations.\\
(a) If $A$ has no complex multiplication then $\NS(A)_{\Q} \cong \ad^0(\cB)$.\\
(b) If $A$ has complex multiplication then $\NS(A)_{\Q} \cong \ad^0(\cB) \oplus \Norm$ where $\Norm$ denotes the one-dimensional representation of $\cB^*$ given by the reduced norm $Nrd: \cB^* \to \Q$. In this case the two summands are orthogonal to each other under the intersection pairing on $\NS(A)_{\Q}$.
\end{lemma}

{\em Proof.} See \cite{besser}. 
\enddemo

We fix now (for the rest of this paper) an imaginary quadratic number field $K$. We denote the non-trivial automorphism of $K$ by $a \mapsto \bar{a}$. We assume that the following conditions hold.

\noi $\bu$ The discriminant $D_K>0$ of $K$ is relatively prime to $N = N^+ N^{-} p$.

\noi $\bu$ All prime divisors of $N^{-}p$ are inert in $K$.

\noi $\bu$ All prime divisors of $N^+$ are split in $K$.

Note that the second condition implies that the prime ideal $p\cO_K$ splits completely in the Hilbert class field $H/K$. For simplicity we also assume that $\cO_K^* = \{ \pm 1 \}$, a condition which is satisfied as soon as $D_K>4$. 

In the case where $A$ has complex multiplication by $\cO_K$ we can pick a canonical generator (up to sign) of the one-dimensional subspace of $\ep_2 \NS(A)_{\Q} \break \seq \NS(A)_{\Q}$ where $\cB^*$ acts via the reduced norm.

\begin{prop} 
\label{prop:prop1} 
Assume that $A$ has complex multiplication by $\cO_K$. Then there exists an element $y_{H}$ in $\NS(A)$ such that\\
(a) $\io(b)^*(y_{CM}) = N(b) y_{CM}\qquad \forall \,\,b\in \cR$\\
(b) The self-intersection number of $y_{CM}$ is $= - 2 D_K$ (thus $y_{CM}\in \ep_2 \NS(A)_{\Q}$).\\
Up to sign, $y_{CM}$ is uniquely determined by these properties. 
\end{prop}

{\em Proof.} The uniqueness of $\pm y_{CM}$ follows immediately from (\ref{lemma:l5}). For the existence we note first that 
\[
\End(A) \cong \cO_K\otimes \cR^{\max} \cong M_2(\cO_K).
\]
Once we have fixed a bijection $\End(A) \stackrel{\simeq}{\lra} M_2(\cO_K)$ we get an isomorphism $A \cong E\times E$ where $E$ is an elliptic curve over $H$ with $\End(E) \cong \cO_K$. Denote by $\G_{\sqrt{-D_K}}$ the graph of the endomorphism $\sqrt{-D_K}\in \End(E)$. Then
\[
y_{CM} \colon = [\G_{\sqrt{-D_K}}] - D_K [0\times E] - [E\times 0]\in \NS(A)
\]
has the desired properties (see \cite{nekovar3}, II.3.3(3)). \enddemo

A Heegner point on $X$ associated to $K$ is a Heegner point which is represented by a triple $(A, \iota, C)$ where $A$ has complex multiplication by the ring of integers $\cO_K$ of $K$ (and $C$ is stable under the $\cO_K$-action). Since from now on we are dealing only with Heegner points associated to $K$ we often drop the reference to $K$.

Let $t$ denote the number of prime factors of $N$ and $h$ the class number of $K$. There are precisely $2^t h$ different Heegner points on $X$ associated to $K$ which are all defined over the Hilbert class field $H$ of $K$ (cf.\ \cite{bertolini_darmon1}, 2.5). Let $\cW\cong (\Z/2\Z)^t$ denote the subgroup of $\Aut(X)$ generated by the Atkin-Lehner involutions $w_{\ell}: X \to X$, for $\ell$ a prime factor of $N$ and let $\Delta= \Pic(\cO_K)$ which we identify with $\Gal(H/K)$ via the Artin map. There is a natural action of $\cW\times \Delta$ on the set of Heegner points which is free and transitive.

To define the Heegner cycle for $\cM_n$ we have to work again with the finer moduli problem \ref{definition:qmlevelM}. We fix an integer $M \geq 3$ relatively prime to $N$ and let $q: X_M \to X$ be the projection and $\pi: \cA \to X_M$ be the universal abelian surface. Let $P$ be a Heegner point on $X$ (viewed as a closed point on $X_H$) and let $\wt{P}$ be any point of $(X_M)_H$ above $P$. The fiber $\cA_{\wt{P}}$ is an abelian surface with $\End_{\cR^max}(\cA_{\wt{P}}) = \cO_K$. Hence by proposition \ref{prop:prop1} there exists an -- up to sign -- unique element $y_{CM, \wt{P}}\in \NS(\cA_{\wt{P}})$ satisfying (a) and (b) of \ref{prop:prop1}. We choose a representative $z_{\wt{P}}$ of $y_{CM, \wt{P}}$ in $\ep_2 \Pic(\cA_{\wt{P}}) = \ep_2 CH^1(\cA_{\wt{P}})$. 

One can choose the elements $z_{\wt{P}}$, $\wt{P}\in q^{-1}(P)$ in such a way that they are compatible with the $G = G_M$-action. More precisely any $g\in G$ extends uniquely to an automorphism $\cA\to \cA$, also denoted by $g$, and thus induces an isomorphism
\[
g: \cA_{\wt{P}}\lra \cA_{g(\wt{P})}
\]
for $\wt{P}\in q^{-1}(P)$. We require that the elements $z_{\wt{P}}\in \ep_2 CH^1(\cA_{\wt{P}})$ satisfy 
\begin{equation}
\label{equation:comp1}
g_*(z_{\wt{P}}) = z_{g(\wt{P})}
\end{equation}
for all $\wt{P}\in q^{-1}(P)$ and $g\in G$.

We define the cycle class $y_P = y_P^{(2)}$ in $CH^2((\cM_2)_H)$ to be the image of $z_{\wt{P}}$ under
\[
\ep_2 CH^1(\cA_{\wt{P}}) \stackrel{(\io_{\wt{P}})_*}{\lra} \ep_2 CH^2(\cA_H) \stackrel{(p_G)_*}{\lra} (\ep_2 CH^2(\cA_H))^{G} = CH^2((\cM_2)_H)
\]
with $p_G = \frac{1}{|G|} \sum_{g\in G} g\in \Q[G]$. By (\ref{equation:comp1}) this is independent of the point $\wt{P}$ we have chosen above $P$.

We also require that the elements $\{y_P\}$ are compatible with respect to the action of the Atkin-Lehner involutions and the $\De$-action. Any $w\in \cW$ extends canonically to an involution $w: X_M \to X_M$ which commutes with the $G$-action. Also the action of $\Delta$ on $H$ induces an action on $(X_M)_H$, $\cA^m_H$ etc.\ which commutes with the $\cW$- and $G$-operations. We can choose the elements $z_{\wt{P}}$ for different Heegner points $P$ and $\wt{P}\in q^{-1}(P)$ so that 
\[
w_*(z_{\wt{P}}) = z_{w(\wt{P})}, \qquad \delta_*(z_{\wt{P}}) = z_{\delta(\wt{P})}
\]
for all $w\in \cW$ and $\delta \in \De$ and all $\wt{P}$. With this normalisation we obtain
\[
w_*(y_P) = y_{w(P)},\qquad \delta_*(y_P) = y_{\delta(P)} \qquad \forall\, w\in \cW, \delta \in \De.
\]

\begin{remark} 
{\rm Albeit the cycle class $y_P\in CH^2((\cM_n)_H)$ is defined using $z_{\wt{P}}\in CH^1(\cA_{\wt{P}})$ conjectures of Bloch and Beilinson would imply that it depends only on the the class of $z_{\wt{P}}$ modulo homological equivalence i.e.\ on $y_{CM, \wt{P}}$.}
\end{remark}

For arbitrary (even) $n\geq 2$ we define cycles $y_P^{(n)}\in CH^{m+1}((\cM_n)_H)$ as follows. Firstly for $\wt{P}\in q^{-1}(P)$ let $z_{\wt{P}}^{(n)}$ be the image of the $m$th exterior product $z_{\wt{P}}\times \ldots \times z_{\wt{P}}$ under the projector $\ep_n$ acting on $CH^m(\cA^m_{\wt{P}})$. Then $y_P^{(n)}\in CH^{m+1}((\cM_n)_H)$ is defined as the image of $z^{(n)}_{\wt{P}}$ under
\[
\ep_n CH^m(\cA^m_{\wt{P}}) \stackrel{(\io_{\wt{P}})_*}{\lra} \ep_n CH^{m+1}(\cA^m_H) \stackrel{(p_G)_*}{\lra} CH^{m+1}((\cM_n)_H)
\]
As before $y_P^{(n)}$ does not depend on the point $\wt{P}$ above $P$ and we have
\begin{equation}
\label{equation:comp2}
w_*(y_P^{(n)}) = y_{w(P)}^{(n)},\qquad \delta_*(y_P^{(n)}) = y_{\delta(P)}^{(n)}
\end{equation}
for all $w\in \cW$, $\delta \in \De$ and all Heegner points $P$.
The cycle class $y^{(n)}_P$ is called Heegner cycle on $\cM_n$.\\

We want to compute the image of $y^{(n)}_P$ under the map $\rho_H$ (cf.\ (\ref{eqn:s7paj})). For that we need to recall the $p$-adic analytic description of the set of Heegner points on $X$ given in (\cite{bertolini_darmon2}, section 5). Let $\cO \colon = \cO_K[\frac{1}{p}]$. An embedding $\Psi: K \to B$ is called optimal if $\Psi(\cO) \seq R$. The group $\G$ acts naturally by conjugation on the set of optimal embeddings and we denote by $\emb(\cO, R)$ the set of conjugacy classes (note that contrary to \cite{bertolini_darmon2} we do not consider orientations here). For an optimal embedding $\Psi$ we denote by $[\Psi]\in \emb(\cO, R)$ its class. There is a natural action of $\Pic(\cO) \cong\De$ on $\emb(\cO, R)$ given as follows (see \cite{bdis}, 2.3). Let $\alpha\in \De$ and choose an ideal $\fa\in \cO$ representing it. Let $\Psi: K \to B$ be an optimal embedding. The right order of the left $R$-ideal $R\Psi(\fa)$, denoted by $R_{\fa}$, is an Eichler $\Z[\frac{1}{p}]$-order of level $N^+$ in $B$. The right action of $\Psi(\cO)$ on $R\Psi(\fa)$ yields an embedding $\wt{\Psi}_{\fa}: \cO \to R_{\fa}$. Since all Eichler $\Z[\frac{1}{p}]$-orders are conjugated there is an element $a\in B^*$ such that $\ord_p(Nrd(a)) = 0$ and $R = a R_{\fa} a^{-1}$. Then, 
\begin{equation}
\label{eqn:actiondelta}
\alpha \cdot [\Psi] \colon = [\wt{\Psi}_{\fa}^a].
\end{equation}

By (\cite{bertolini_darmon2}, section 5) there is a bijection
\begin{equation}
\label{eqn:pheeg}
\mbox{Heegner points on $X_H$} \,\, \stackrel{\cong}{\longleftrightarrow} \emb(\cO,R)
\end{equation}
which is compatible with the $\De$-action on both sides. It is given as follows. We fix an embedding $\io: H(\mu_M) \hookrightarrow \Qur$ and denote by $v$ be the corresponding place of $H(\mu_M)$ (and also of $H$). The embedding $\io$ allows us to identify $K_p \colon = K\otimes \Q_p$ and $H_v$ with $\Q_{p^2} \seq \Qur$. Using the $p$-adic uniformisation \ref{theorem:cer_drinf} of $X_{\Q_{p^2}}$ we view $\emb(\cO,R)$ as a subset of $X(K_p)$ via 
\begin{equation}
\label{eqn:pheeg2}
\emb(\cO,R) \hookrightarrow \G\backslash\Hom(K_p, B_p) \cong \G\backslash\Hom(\Q_{p^2}, M_2(\Q_p)) = X_{\G}(K_p)
\end{equation}
The image of $\emb(\cO,R)$ under (\ref{eqn:pheeg2}) is the set of Heegner points (considered as a subset of $X(K_p)$ via $X(H) \hookrightarrow X(H_v) = X(K_p)$).

Let $P$ be a Heegner point and let $\Psi: K \to B$ be an optimal embedding corresponding to $P$ which we view as an element of $\Hom(\Q_{p^2}, M_2(\Q_p))\seq \cH_p$. We fix one of the components $\G_M\backslash\cH_p$ of $(X_M)_{\Qur}$ and let $\wt{P}\in q^{-1}(P)$ be the point $\G_M z$. Let $P_{\Psi}(X) \in \cP_2$ be the polynomial $P_{\Psi(\sqrt{-D_K})}(X)$ defined as in equation (\ref{eqn:defpu}) in section \ref{sec:dieumod} and define $cl_{\Psi}^{(n)}\in V_n$ by 
\[
cl_{\Psi}^{(n)}(P(X)) \colon = <P(X), P_{\Psi}^m>.
\]

\begin{lemma}
\label{lemma:cyclemapheeg}
Let
\begin{equation}
\label{eqn:s7cycle}
cl_{p, \Psi}: CH^m((\cA^m_{\wt{P}}, \ep_n)) \lra \G((V_n)_{\Psi}[-m])
\end{equation}
denote the upper horizontal homomorphism in diagram (\ref{eqn:cdaj1}). We have, \[
cl_p(z^{(n)}_{\wt{P}}) = \pm cl_{\Psi}^{(n)}.
\]
\end{lemma}

{\em Proof.} It is enough to consider the case $n=2$ since the case of an arbitrary (even) $n$ follows formally from the first case. We have an isomorphism
\begin{equation}
\label{eqn:s7iso1}
\cU_{\Qur} = \{U \in M_2(\Qur)\mid \, \trace(U) = 0\} \lra V_{\Qur}, U \mapsto < P_U(X), \ldots>.
\end{equation}
The Frobenius on $V_{\Psi}[-1]$ corresponds under (\ref{eqn:s7iso1}) to the map $1_{\cU}\otimes \s$ on $\cU_{\Qur}$ and the subspace $V_{\Qur,i}$ to
\[
\cU_i \colon = \{U\in \cU_{\Qur}\mid\, \Psi(\bar{a}) U \Psi(a)  = a^i \,\bar{a}^{2-i} U\}.
\]
A simple computation now shows that $\G(V_{\Psi}[-1])$ is an one-dimensional $\Q_p$-vector space generated by $<P_{\Psi}, \ldots> = cl_{\Psi}^{(2)}$. Since $< P_{\Psi}, P_{\Psi}> = -\trace( \sqrt{-D_K} \overline{(\sqrt{-D_K})}) = -2 D_K$ the assertion follows from \ref{remark:dieumodformpair2} and \ref{prop:prop1}. \enddemo

\begin{remark} 
\label{remark:dieumodformsign}
{\rm Recall that the isomorphism (\ref{eqn:dieumodformn2}) was so far well-defined only up to sign (Remark \ref{remark:dieumodformpair}) and that changing the sign also changes the sign of the map (\ref{eqn:s7cycle}). We fix the sign by requiring that $cl_p(z^{(n)}_{\wt{P}}) = cl_{\Psi}^{(n)}$. This automatically implies that $cl_p(z^{(n)}_{\wt{P}}) = cl_{\Psi}^{(n)}$ for {\it any} optimal embedding $\Psi:K \to B$ since $\cW\times \De$-acts transitively on the set of Heegner points. For example consider the Atkin-Lehner involution $w_p$ at $p$ and choose $\gamma\in R^*$ such that $\ord_p(Nrd(\gamma))$ is odd. We have a commutative diagram
\[
\begin{CD}
CH^m((\cA^m_{\wt{P}}, \ep_n)) @> cl_{p, \Psi} >> \G((V_n)_{\Psi}[-m])\\
@VV (w_p)_* V @VV \gamma_* V\\
CH^m((\cA^m_{w_p(\wt{P})})), \ep_n)) @> cl_{p, \Psi^{\gamma}} >> \G((V_n)_{\Psi^{\gamma}}[-m])
\end{CD}
\]
where $\gamma_*(v) \colon = Nrd(\gamma)^{-m} \gamma v$. Since $\gamma_*(cl_{\Psi}^{(n)}) = cl_{\Psi^{\gamma}}^{(n)}$ we get 
\[
cl_{p, \Psi^{\gamma}}(y^{(n)}_{\wt{w_p(P)}}) = (w_p)_*(y^{(n)}_{\wt{P}}) = \gamma_*(cl_{p, \Psi}^{(n)}) = cl_{\Psi^{\gamma}}^{(n)}.
\]
For $\alpha \in \Delta$ and $a\in B^*$ as in (\ref{eqn:actiondelta}) we get a diagram as above with $\alpha$ and $a$ instead of $w_p$ and $\gamma$.}
\end{remark}

\section{Derivatives of anticyclotomic $p$-adic $L$-func\-tions}
\label{sec:lfunc}

In this section we show that the derivative of the anticyclotomic $p$-adic $L$-function attached to a newform of even weight $k \geq 4$ on $X$ at the central critical point $s=k/2$ can be expressed in terms of the $p$-adic Abel-Jacobi image of Heegner cycles. 

To begin with we briefly review the construction of the anticyclotomic $p$-adic $L$-function $L_p(f, K, s)$ from (\cite{bdis} section 2.5). We keep the notation and assumptions of the last section. Let $f(z) \in M_k(\G)$ and let $\cA_n$ denote the set of $\C_p$-valued functions on $\PP^1(\Q_p)$ which are locally analytic except for a pole of order at most $n$ at $\infty$. Associated to $f(z)$ there is a distribution
\[
\mu_f: \cA_n \to \C_p, \vartheta \mapsto \mu_f(\vartheta) = \int_{\PP^1(\Q_p)} \, \vartheta(x) \mu_f(dx)
\]
which is characterised by the property
\[
\int_{U(e)}\, P(x) \mu_f(dx) = \Res_e(P(z) f(z) dz)
\]
for all $P(X) \in \cP_n$ and $e\in \edges(\cT)$. Here $U(e)$ denotes as usual the compact subset of $\PP^1(\Q_p)$ corresponding to the {\it ends} of $\cT$ containing $e$ (see e.g.\ \cite{bdis}, section 1.3). 

Let $K_{\infty}/K$ be the maximal abelian extension of $K$ which is unramified outside $p$ and {\it anticyclotomic}, i.e.\ the involution in $\Gal(K/\Q)$ acts as $-1$ on $\Gal(K_{\infty}/K)$. Let $G\colon = \Gal(K_{\infty}/H)$ and $\wt{G} \colon = \Gal(K_{\infty}/K)$ so that $\wt{G}/G \cong \De$. Class field theory identifies $G$ with $K_p^*/\Q_p^*$. An optimal embedding $\Psi: K \to B$ yields an embedding
\[
G \cong K_p^*/\Q_p^* \stackrel{\Psi}{\lra} B_p^*/\Q_p^* \cong \PGL_2(\Q_p)
\]
hence a (simply transitive) action of $G$ on $\PP^1(\Q_p)$. For a fixed base point $\star\in \PP^1(\Q_p)$ we let $\eta_{\Psi, \star}: G\to\PP^1(\Q_p)$ be the bijection given by letting $G$ act on $\star$. We obtain a distribution $\mu_{f,\Psi,\star}$ on $G$, i.e.\ a functional on the set $\cA(G)$ of $\C_p$-valued locally analytic functions on $G$ by setting
\[
\int_{G} \, \vartheta(\alpha) \mu_{f,\Psi,\star}(d\alpha) \colon = \int_{\PP^1(\Q_p)} \, \vartheta\circ \eta_{\Psi, \star}^{-1}(x) P_{\Psi}(x)^m \mu_f(dx)
\]
for $\vartheta\in \cA(G)$. For $\ga\in \G$ we have $\mu_{f,\Psi^{\ga},\ga\star} =\mu_{f,\Psi,\star}$. The distribution $\mu_{f,\Psi,\star}$ depends on the base point $\star$ and the representative $\Psi$ of $[\Psi]\in \emb(\cO, R)$ only up to translation by an element in $G$. 

Define $\log: G \to \Qur$ to be the composite
\begin{equation}
\label{eqn:log}
\log: G \cong K_p^*/\Q_p^* \cong K_{p,1} \stackrel{\log_p}{\lra} \Q_{p^2} \seq \Qur
\end{equation}
where $K_{p,1} = \{x\in K_p^*\mid\, N_{K_p/\Q_p}(x) = 1\}$ and where the second isomorphism is given by $x\mapsto x/\bar{x}$. The map (\ref{eqn:log}) extends uniquely to $\wt{G}$. For $\alpha$ in $G$ or $\wt{G}$ and $s\in \Z_p$ we define $\alpha^s \colon = \exp(s \,\log(\alpha))$. The partial $p$-adic $L$-function attached to the datum $(f, \Psi, \star)$ is defined as 
\[
L_p(f, \Psi, \star, s) \colon = \int_{G} \, \alpha^{s-k/2} \mu_{f,\Psi,\star}(d\alpha)
\]
If $\star'\in\PP^1(\Q_p)$ is another base point and $\Psi': K \to B$ another optimal embedding which is conjugate under $\G$ to $\Psi$ then there exists an element $\alpha_0\in G$ such that
\[
L_p(f, \Psi', \star', s) = \alpha_0^{s-k/2} L_p(f, \Psi, \star, s).
\]
In particular since $L_p(f, \Psi, \star, k/2) = 0$ we see that the first derivative $L_p'(f, \Psi,\break \star, k/2)$ is independent of the choice of $\Psi$ in $[\Psi]$ and $\star$.

The distribution $\mu_{f,\Psi, \star}$ extends canonically to a distribution on $\wt{G}$ denoted by $\mu_{f,K}$ for short. For that let $\La$ denote the set of $\G$-conjugacy classes of pairs $(\Psi,\star)$ where $\Psi: K \to B$ is an optimal embedding and $\star\in \PP^1(\Q_p)$. The action of $\cW\times \Delta$ on $\emb(\cO, R)$ lifts canonically to a simply transitive action of $\cW\times\wt{G}$ on $\La$ such that $\alpha (\Psi, \star) = (\Psi, \Psi(\alpha) \star)$ for $\alpha\in G$ and $(\Psi, \star) \in \La$ (see \cite{bdis}, Lemma 2.14). We denote by $\Xi$ the set of $\De$-orbits in $\emb(\cO, R)$ (or $\wt{G}$-orbits in $\La$). It is a principal homogeneous $\cW$-space. 

Fix $\xi$ in $\Xi$ and a pair $(\Psi,\star)$ representing it (depending on the circumstances we consider elements in $\Xi$ either as $\De$-orbits in $\emb(\cO, R)$ or $\wt{G}$-orbits in $\La$). For every $\delta\in\De$ we choose a lift $\alpha_{\delta}\in \wt{G}$ and write $(\Psi_{\delta},\star_{\delta}) \colon = \alpha_{\delta}(\Psi,\star)$. The distribution $\mu_{f,K} = \mu_{f,K,\xi}$ is given by the formula
\[
\int_{\wt{G}}\, \vartheta(\alpha) \mu_{f,K}(d \alpha) \colon = \sum_{\delta\in \De} \int_G \,\vartheta(\alpha \alpha_{\delta}) \mu_{f,\Psi_{\delta},\star_{\delta}}(d\alpha)
\]
for $\vartheta\in \cA(\wt{G})$. It is independent of the choice of the $\alpha_{\delta}$ and depends on the pair $(\Psi,\star)\in \xi$ only up to translation by an element in $\wt{G}$. The anticyclotomic $p$-adic $L$-function attached to the datum $(f, K, \xi)$ is a function in the $p$-adic variable $s\in \Z_p$ defined as 
\[
L_p(f, K, \xi, s) \colon = \int_{\wt{G}} \, \alpha^{s-k/2} \mu_{f,\Psi,\star}(d\alpha).
\]
Up to a multiple of the form $s\mapsto \alpha_0^{s-k/2}$ it is independent of the pair $(\Psi,\star)\in\xi$ (for the justification of the term anticyclotomic $p$-adic $L$-function we refer to \cite{bdis}). $L_p(f, K, \xi, s)$ can be written as a sum of partial $p$-adic $L$-functions
\begin{equation}
\label{eqn:sumpartial}
L_p(f, K, \xi, s) = \sum_{i=1}^h \, L_p(f, \Psi_i, \star_i, s)
\end{equation}
where $\Psi_1, \ldots, \Psi_h$ are representatives of the $\De$-orbit $\xi$ and $\star_1, \ldots, \star_h$ are suitably chosen base points. In particular $L_p(f, K, \xi, k/2) = 0$ and $L_p'(f, K, \xi, k/2) \break = \sum_{i=1}^h \, L_p'(f, \Psi_i, \star_i, k/2)$ is independent of the pair $(\Psi,\star)$.

\begin{remark} 
\label{remark:xi} {\rm Assume that $f(z)$ is a newform, i.e.\ it corresponds to a newform $f_{\infty}\in S_k(\G_0(N))$ under the isomorphism (\ref{eqn:jacquet}). In this case $L_p(f, K, \xi, s)$ is independent of $\xi$ up to sign. In fact if $w\in \cW$ and $\ep\in \{\pm 1\}$ such that $w(f) = \ep f$ then 
\[
L_p(f, K, w \xi, s) = \ep L_p(f, K, \xi, s).
\]
In \cite{bdis} the dependence of $L_p(f, K, \xi, s)$ on $\xi$ was surpressed. We had considered there always a fixed $\De$-orbit $\xi$ in $\emb(\cO, R)$. However the reader should notice that there is no canonical choice for it.}
\end{remark}

For $[\Psi]\in \emb(\cO, R)$ we define
\[
\cL_p'(\Psi, k/2): M_k(\G) \lra \Qur, f(z) \mapsto L_p'(f, \Psi, \star, k/2)
\]
and for $\xi \in \Xi$ we let $\cL_p'(K,\xi,k/2)$ be the functional
\[
M_k(\G) \lra \Qur, f(z) \mapsto L_p'(f, K, \xi, k/2),
\]
so that $\cL_p'(K, \xi, k/2) = \cL'(\Psi_1, k/2) + \ldots + \cL'(\Psi_h, k/2)$ where $\Psi_1, \ldots, \Psi_h$ are representatives of $\xi$. Let $P_1,\ldots, P_h$ be the Heegner points on $X_H$ corresponding to $\Psi_1, \ldots, \Psi_h$ via (\ref{eqn:pheeg}) and let $\Pbar_i\in X(H)$ be the complex conjugate of $P_i$. We denote by $y^{(n)}_{\xi}$ the image of the Heegner cycle $y_{P_i}^{(n)}$ under the push-forward $CH^{m+1}((\cM_n)_H) \to CH^{m+1}((\cM_n)_K)$. It does not depend on the choice of $i\in \{1, \ldots, h\}$ and is mapped to $y^{(n)}_{P_1}  + \ldots + y^{(n)}_{P_h}$ under pull-back to $H$. Complex conjugation also acts on the set of optimal embeddings (by $\Psibar(a) \colon = \Psi(\bar{a})$) and likewise on $\emb(\cO, R)$ and $\Xi$.

Our main result is:

\begin{theorem}
\label{theorem:main}
\[
\cL_p'(K, \xi, k/2) = \rho_H(y^{(n)}_{P_1}  + \ldots + y^{(n)}_{P_h} + (-1)^{m+1} (w_p)_*(y^{(n)}_{\Pbar_1}  + \ldots + y^{(n)}_{\Pbar_h})) = 
\]
\[
=\rho_{K}(y^{(n)}_{\xi}  + (-1)^{m+1} (w_p)_*(y^{(n)}_{\bar{\xi}})).
\]
\end{theorem}

{\em Proof.} The second equality follows immediately from the definition of $y^{(n)}_{\xi}$. For the first equality it is enough to show that
\begin{equation}
\label{eqn:main}
\cL_p'(\Psi, k/2) = \rho_H(y^{(n)}_P + (-1)^{m+1} y^{(n)}_{w_p(\Pbar)})
\end{equation}
for any Heegner point $P$ and corresponding optimal embedding $\Psi: K \to B$.
Fix $P, \Psi$ and let $z_0$ and $\bar{z}_0$ be the points on $\cH_p$ which correspond to $\Psi$ and $\Psibar$ respectively via (\ref{eqn:puh}). Note that the class $[\Psibar]\in \emb(\cO, R)$ belongs to $w_p(\Pbar)$ via (\ref{eqn:pheeg}) (see \cite{drinfeld, boutot_carayol}). Let $U \colon = X-\{P, w_p(\Pbar)\}$ and let 
\begin{equation}
\label{eqn:s8ext}
0 \lra H^1(X_{\G}, \cV_n) \lra E \lra \Qur[m+1] \lra 0
\end{equation} 
be the pull-back of 
\[
0 \lra H^1(X_{\G}, \cV_n) \stackrel{j_*}{\lra} H^1(U_{\Qur}, \cV_n) \lra (V_n)_{\Psi}[1] \oplus (V_n)_{\Psibar}[1] \lra 0
\]
with respect to 
\[
\Qur[m] \lra  (V_n)_{\Psi} \oplus (V_n)_{\Psibar}, 1 \mapsto (cl_{\Psi}^{(n)}, - cl_{\Psibar}^{(n)}).
\]
By Proposition \ref{prop:ajfiber}, Lemma \ref{eqn:s7cycle} and Remark \ref{remark:dieumodformsign} the extension class (\ref{eqn:s8ext}) corresponds to $\rho_H(y^{(n)}_P + (-1)^{m+1} y^{(n)}_{w_p(\Pbar)})$ under the isomorphism
\[
\Ext^1_{MF^{ad}_{\Qur}(\p, N)}(\Qur[m+1], H^1(X_{\G}, \cV_n)) \cong H^1(X_{\G}, \cV_n)/F^{m+1}\cong M_k(\G)^{\vee}.
\]

Concretely, let $\alpha$ be the uniquely determined element in $H^1(U_{\Qur}, \cV_n)_{m+1}$ with $N(\alpha) =0$ and such that
\[
\Res_{z_0}(\alpha) = cl_{\Psi}^{(n)} = < P_{\Psi}^m, \ldots>,\quad \Res_{\bar{z}_0}(\alpha) = (-1)^{m+1} cl_{\Psibar}^{(n)} = -< P_{\Psi}^m, \ldots>.
\]
Let $\beta\in H^1(X_{\G}, \cV_n)$ with $j_*(\beta) \equiv \alpha \,\mod{F^{m+1}}$. Then,
\begin{equation}
\label{eqn:s8aj1}
\rho(y^{(n)}_P + (-1)^{m+1} y^{(n)}_{w_p(\Pbar)})(f(z)) = < [\omega_f], \beta>
\end{equation}
for all $f(z) \in M_k(\G)$. 

Since $H^1(X_{\G}, \cV_n) = H^1(X_{\G}, \cV_n)_m \oplus F^{m+1}H^1(X_{\G}, \cV_n)$ (see the proof of Theorem \ref{theorem:lfm=lt}) we may assume that $\beta\in H^1(X_{\G}, \cV_n)_m = \io(H^1(\G, (V_n)_{\Qur}))$, i.e.\ $\beta = \io(c)$ for some $c\in H^1(\G, (V_n)_{\Qur})$. In order to compute $< \beta, [\omega_f]>$ we use the results and notation of appendix \ref{subsection:a2}. By Theorem \ref{theorem:deshalit} we have 
\[
< [\omega_f], \beta> = - <I([\omega_f], c>_{\G}.
\]
Let $\chi$ be a $\G$-invariant $V_n$-valued meromorphic differential form on $\cH_p$ which is holomorphic outside of $\pi^{-1}(U^{an})$, which has simple poles at $z_0, \bar{z}_0$ and whose class in $F^{m+1} H^1(U_{\Qur}, \cV_n)$ represents $\alpha - j_*(\beta)$. Hence $P_U([\chi]) = - P_U(j_*(\beta)) = -c, N([\chi]) = 0$ and
\[
\Res_{z_0}(\chi) = \Res_{z_0}(\alpha) = < P_{\Psi}^m, \ldots> = - \Res_{\bar{z}_0}(\alpha) = - \Res_{\bar{z}_0}(\chi).
\]
We can apply Corollary \ref{coro:openshalit} and obtain
\[
< [\omega_f], \beta> = - <I([\omega_f]), c>_{\G} = <I([\omega_f]),P_U([\chi])>_{\G} =
\]
\[
= < F_{\omega_f}(z_0), \Res_{z_0}(\chi)> + 
< F_{\omega_f}(\bar{z}_0), \Res_{\bar{z}_0}(\chi)> = \int^{z_0}_{\bar{z}_0}\,\, f(z) P_{\Psi}(z)^m dz.
\]
By (\cite{bdis}, Theorem 3.5) the last expression is equal to $L_p'(f, \Psi, \star, k/2)$ hence together with (\ref{eqn:s8aj1}) we get
\[
\rho(y^{(n)}_P + (-1)^{m+1} y^{(n)}_{w_p(\Pbar)})(f(z)) = L_p'(f, \Psi, \star, k/2)
\] 
as claimed.\enddemo

Let $f(z)\in M_k(\G, \C_p)$ be a newform and $w$ be the eigenvalue of the Atkin-Lehner involution $w_p$ acting on $f(z)$ (so $w= -1$ if $f(z)$ is of split multiplicative type and $w=1$ otherwise). Let
\[
\rho_{f,p}:  CH^{m+1}((\cM_n)_K)\lra \C_p, z \mapsto \rho_K(z)(f(z))
\]
be the $f(z)$-component of $\rho_K$.

\begin{coro}
\label{coro:main}
We have the following equality
\[
L_p'(f, K,\xi, k/2) = \rho_{f,p}(y^{(n)}_{\xi}  + w (-1)^{m+1} y^{(n)}_{\bar{\xi}}).
\]
\end{coro}

\begin{remark} 
\label{remark:selmer}
{\rm Let $f(z)\in M_k(\G, \C_p)$ be a newform. In \cite{bloch_kato} the authors have defined a subgroup $H^1_f(K, V_p(f))$ of $H^1_{\cont}(K, V_p(f))$ which one might consider as a cohomological version of a Selmer group. Corollary \ref{coro:main} together with Kolyvagin's method of Euler systems (compare \cite{nekovar1}) can be used to show that 
\[
L_p'(f, K,\xi, k/2) \neq 0 \Longrightarrow \dim(\Image(cl_{0,f})) = \dim(H^1_f(K, V_p(f))) = 1.
\]
Here $cl_{0,f}: CH^{m+1}((\cM_n)_K \to H^1_{\cont}(K, V_p(f))$ is the $f$-component of the cohomological Abel-Jacobi map (\ref{eqn:aj1}) (for the motive $\cM_n$). This statement can be interpreted as an affirmative answer to a Bloch-Beilinson type conjecture relating the order of vanishing of $L_p(f, K,\xi, s)$ at $s=k/2$ to the dimension of the Abel-Jacobi image of a cycle class group attached to the motive ``$M(f)$''. The reader is invited to formulate a precise conjecture (generalising conjecture 4.1 of \cite{bertolini_darmon1}).}
\end{remark}

\section{Appendix}

\subsection{Chow motives attached to modular forms on Shi\-mura curves}
\label{subsection:a1}

Attached to the space of modular forms of weight $k = n+2$ on $X$ is a certain motive $\cM_n$ which is the analogue in the context of Shimura curves of the motive attached to the space of cusp forms of a fixed weight and level constructed by Scholl \cite{scholl}. The motive $\cM_n$ has been considered first by Besser (see \cite{besser}). He construced it as a motive for absolute Hodge cycles. It is possible to refine his construction and define $\cM_n$ in the category of Chow motives (see \cite{wortmann}). For convenience we briefly recall here the main steps of the construction and some properties of $\cM_n$.

Let $K$ be a field of characteristic 0 and $S$ be a smooth quasiprojective connected variety over $K$. We denote by $\cM^0(S)$ the category of relative Chow motives over $S$ with respect to graded correspondences (see \cite{de-mu, kuennemann}). Recall that the objects of $\cM^0(S)$ are triples $(X, p, i)$ where $X$ is a smooth projective $S$-scheme, $p$ is a projector (i.e.\ an idempotent in the ring $\Corr_S(X, X) \colon = \bigoplus_{\nu} \, CH^{\dim(X_{\nu}/S)}(X_{\nu}\times_S X)$ of relative correspondences where $X = \coprod_{\nu} \, X_{\nu}$ is the decomposition of $X$ into connected components) and $i$ is an integer. Morphisms are given by
\[
\Hom_{\cM^0(S)}((X, p, i), (Y, q, j)) = q \circ CH^{\dim(X_{\nu}/S)+j-i}(X_{\nu}\times_S Y) \circ p.
\]
Composition is induced by composition of correspondences. The category $\cM^0(S)$ is an additive pseudo-abelian $\Q$-linear rigid tensor category. In particular the kernel of projector exists in $\cM^0(S)$. For a smooth projective $S$-scheme $X$ we denote by $h(X/S)\colon = (X, \De_X)$ its motive (we abbreviate $(X, p)$ for $(X, p, 0)$). 

Let $S = X_M/\Q$ be the fine moduli scheme of the moduli problem \ref{definition:qmlevelM} for some $M \geq 3$ and let $\pi:\cA \to X_M$ be the universal abelian surface with quaternionic multiplication. The motive $h(\cA) = h(\cA/X_M)$ admits a canonical decomposition
\[
h(\cA) = h^0(\cA) \oplus h^1(\cA) \oplus h^2(\cA)\oplus h^3(\cA)\oplus h^4(\cA)
\]
with $h^i(\cA) \cong \Lambda^i h^1(\cA)$ and $h^i(\cA)^{\vee} \cong h^{4-i}(\cA)(2)$ (see \cite{de-mu} or \cite{kuennemann}). The $\cR^{\max}$-action induces an embedding $\cB \hookrightarrow \End(h^1(\cA))$. There exists a unique idempotent $e\neq 0$ in the subring $\sym^2 \cB$ of $\End(\Lambda^2 h^1(\cA))$ with $x^2 \cdot e = Nrd(x) e$ for all $x\in \cB$. The projector $\ep_2 \in \Corr_{X_M}(\cA, \cA)$ is defined by $(\cA, \ep) = h^2(\cA)_{-} \colon = \Ker(e)$.

Let $<\,\,\,, \,\,\,>: \sym^2 h^2(\cA)_{-} \to \Q(2)$ be the restriction of the symmetric pairing 
\[
h^2(\cA) \otimes h^2(\cA)\cong \Lambda^2 h^1(\cA) \otimes\Lambda^2 h^1(\cA) \stackrel{\wedge}{\lra}\Lambda^4 h^1(\cA) \cong \Q(-2)
\]
to $h^2(\cA)_{-}$. Let 
\[
\De_m: \sym^m h^2(\cA)_{-} \lra (\sym^{m-2} h^2(\cA)_{-})(-2)
\] 
be the Laplace operator associated to $<\,\,\,, \,\,\,>$ and let $\lambda_{m-2}: \sym^{m-2} h^2(\cA)_{-}\to (\sym^m h^2(\cA)_{-})(2)$ be the map given symbolically by $\lambda_{m-2}(x_1 x_2 \ldots x_{m-2}) = x_1 x_2 \ldots x_{m-2} \mu$ where $\mu: \Q\to (\sym^2 h^2(\cA)_{-})(2)$ is the dual of $<\,\,\,, \,\,\,>$ (twisted by $2$). Then $\De_m \circ\la_{m-2}$ is an isomorphism. Put $p \colon = \la_{m-2}\circ(\De_m \circ\la_{m-2})^{-1}\circ\De_m$. Since $p^2 = p$, $\Ker(p)$ exists (hence also the kernel of $\De_m$ exist and is equal to $\Ker(p)$). The projector $\ep_n\in \Corr_{X_M}(\cA^m, \cA^m)$ is defined by $(\cA^m, \ep_n) = \Ker(p)$. 

The $p$-adic realisation of $(\cA^m, \ep_n)$ is equal to the sheaf $\bL_n$ (shifted by $-n$), i.e.\ we have
\begin{equation}
\label{eqn:preal}
R_p(\cA^m, \ep_n) = \bL_n[-n]
\end{equation}
where $R_p: \cM^0(X_M) \to D^b(X_M, \Q_p)$ is the $p$-adic realization functor from $\cM^0(X_M)$ to the bounded derived category of $\Q_p$-sheaves (see \cite{de-mu}, 1.8).

The group $G \colon = (\cR^{\max}/M\cR^{\max})^* \cong \GL_2(\Z/M\Z)$ acts canonically as $X$-automorphisms on $X_M$ and $\cA^m$. We can therefore view the idempotent $p_G \colon = \frac{1}{|G|} \sum_{g\in G} g\in \Q[G]$ as an element in $\Corr_{X}(\cA^m, \cA^m)$. Since the projectors $\ep_n$ and $p_G$ commute their product is a projector, too. The relative Chow motive $\cM_n$ over $X$ is defined as $(\cM_n, p_G \circ \ep_n)$. It is independent of the integer $M$. When considered as a Chow motive over $\Q$ (by applying the canonical functor $\cM^0(X) \to \cM^0(\Q)$) its $p$-adic realization is given by the following lemma.

\begin{lemma}
\label{lemma:a1hpm} We have,
\[
H_p^i(\cM_n) \cong \left\{ \begin{array}{ll}
                  H^1(\overline{X}_M, \bL_n)^{G} & \                  \mbox{if $i=n+1$},\\
                  0 & \mbox{otherwise}.
                \end{array}
\right.
\]
\end{lemma}

To see this note that by (\ref{eqn:preal}) we have
\[
H_p^i(\cM_n) = (p_G)_*(H^i(\overline{X}_M, R_p((\cA^m, \ep_n)/X_M))) = H^{i-n}(\overline{X}_M,\bL_n)^{G}
\]
Since $H^i(\overline{X}_M,\bL_n)=0$ for all $i\neq 1$ the lemma follows. 

Similarly for the de Rham realization we have $H_{\DR}^i(\cM_n) =0$ for $i\neq n+1$ and $F^{m+1}H_{\DR}^{n+1}(\cM_n) \cong M_k(X)$ (see also \ref{proposition:filhodge}).

\subsection{A formula for the cup product on open subschemes of Mumford curves}
\label{subsection:a2}

Let $k$ be a finite extension of $\Q_p$ and $\G\seq \SL_2(k)/\{\pm 1\}$ be a discrete cocompact subgroup. Let $X_{\G}$ be the associated Mumford curve; as a rigid analytic space it is $\G\backslash\cH_p$ where $\cH_p= \PP^1_{/k}-\PP^1(k)$ is the $p$-adic upper half plane defined over $k$. Let $V$ be a finite-dimensional $\G$-representation (over $k$) which carries a $\G$-invariant bilinearform
\begin{equation}
\label{eqn:apair1}
< \,\,\,,\,\,\,> : V \otimes V \lra k
\end{equation}
The locally free $\cO_{\cH_p}$-module $\cO_{\cH_p}\otimes V$ descends to a local system on $X_{\G}$ denoted by $\cV$. The pairing (\ref{eqn:apair1}) induces a pairing
\begin{equation}
\label{eqn:apair2}
< \,\,\,,\,\,\,>_{X_{\G}} : \Hx \otimes \Hx \stackrel{\cup}{\lra} H^2_{\DR}(X_{\G}, \cV\otimes \cV) \lra
\end{equation}
\[
\lra H^2_{\DR}(X_{\G}) \cong k.
\]
In \cite{shalit2} de Shalit proves a formula for (\ref{eqn:apair2}) in terms of ``Coleman'' and ``Schneider integration'' (see Theorem \ref{theorem:deshalit} below). In this section we generalise this formula to certain open subschemes of $X_{\G}$ in the case where $X_{\G}$ is of arithmetic origin (i.e.\ a Shimura curve). Under this hypothesis we give also a new short proof of \ref{theorem:deshalit} based on Theorem \ref{theorem:dieumodform}.

To begin with we recall some facts from \cite{shalit1, shalit3}. There are pairings
\begin{equation}
\label{eqn:apair3}
< \,\,\,,\,\,\,>_{\G} : H^1(\G, V)\otimes C_{\har}(V)^{\G} \lra k,
\end{equation}
\begin{equation}
\label{eqn:apair4}
< \,\,\,,\,\,\,>_{\G} : C_{\har}(V)^{\G}\otimes H^1(\G, V)\lra k
\end{equation}
defined as 
\[
H^1(\G, V)\otimes C_{\har}(V)^{\G}\stackrel{\cup}{\lra} H^1(\G, V\otimes C_{\har}(V))\stackrel{*}{\lra} H^1(\G, C_{\har}(k))\stackrel{tr}{\lra} k,
\]
\[
C_{\har}(V)^{\G}\otimes H^1(\G, V)\stackrel{\cup}{\lra} H^1(\G, C_{\har}(V)\otimes V)\stackrel{*}{\lra} H^1(\G, C_{\har}(k))\stackrel{tr}{\lra} k,
\]
where (*) is induced by (\ref{eqn:apair1}). The map $tr$ can be defined as follows. We have short exact sequences
\begin{equation}
\label{eqn:aseq1}
0 \lra H^i(\G, V) \stackrel{\io}{\lra} H^i_{\DR}(X_{\G}, \cV)\lra H^{i-1}(\G,C_{\har}(V)) \lra 0.
\end{equation}
For $i=2$ the first group vanishes since $\G$ has a subgroup of finite index which is free. Hence for $V=k$ we have an isomorphism 
\[
tr: H^1(\G,C_{\har}(k)) \cong H^2_{\DR}(X_{\G}) \cong k.
\]
Explicitely the pairing (\ref{eqn:apair3}) is given by
\begin{equation}
\label{eqn:apair5}
<[\fz], f>_{\G}\, = \frac{1}{[\G:\G']} \sum_{i=1}^g \, <\fz(\ga_i), f(c_i)>
\end{equation}
for $f\in C_{\har}(V)^{\G}, \fz\in Z^1(\G,V)$. Here $\G'$ is a free subgroup of finite index in $\G$, $b_1 \ldots, b_g, c_1, \ldots, c_g$ are the ``free edges'' of a good fundamental domain $\frak F$ for $\G'\backslash\cT$ (in the sense of \cite{shalit1}, 2.5) and $\ga_1, \ldots, \ga_g$ are generators of $\G'$ with $\ga_i(b_i) = \bar{c}_i$.

\begin{lemma}
\label{lemma:adj} Let $I: \Hx \to C_{\har}(V)^{\G}$ be the map (\ref{eqn:isch}).\\
(a) $< \io(x), y >_{X_{\G}} = < x , I(y) >_{\G}$ for all $x\in H^1(\G,V)$ and $y\in \Hx$.\\
(b) $< y, \io(x) >_{X_{\G}} = - < I(y), x >_{\G}$ for all $x\in \Hx$ and $y\in H^1(\G,V)$.
\end{lemma}

{\em Proof.} We have $H^i_{\DR}(X_{\G}, \cV) = H^i(\G, \Omega^{\bu}(\cH_p)\otimes V)$ and there is a distinguished triangle 
\[
V \stackrel{\io}{\lra} \Omega^{\bu}(\cH_p)\otimes V \stackrel{I}{\lra} C_{\har}(V)[-1] \lra V[1]
\] 
which induces the sequences (\ref{eqn:aseq1}). Part (a) follows from the commutativity of the diagram
\[
\begin{CD}
\Omega^{\bu}(\cH_p)\otimes V \otimes V @>>> \Omega^{\bu}(\cH_p)\\
@VV I \otimes {1_V} V @VV I V\\
C_{\har}(V)[-1]\otimes V @>>> C_{\har}(k)[-1].
\end{CD}
\]
and (b) by applying (a) to the pairing $<v_1,v_2>' \colon = <v_2, v_1>$ instead of (\ref{eqn:apair1}).\enddemo

In \cite{shalit1, shalit2} the following formula is proved through some complicated explicite computations.

\begin{theorem}
\label{theorem:deshalit}
For all $x,y \in \Hx$ we have 
\begin{equation}
\label{eqn:adeshalit}
<x,y>_{X_{\G}} = <P(x),I(y)>_{\G} - <I(x),P(y)>_{\G}.
\end{equation}
\end{theorem}

Before we explain how this can be generalised to open subschemes we show how to deduce \ref{theorem:deshalit} from \ref{theorem:dieumodform} in the case where $X_{\G}$ is the Shimura curve $X$ (over $\Qur$), $V= V_n$ with $n =2m \geq 2$ even and (\ref{eqn:apair1}) is the pairing (\ref{eqn:s4pair4}).

Under these assumptions the pairing (\ref{eqn:apair2}) can be interpreted as a map in $MF_{\Qur}(\p, N)$ (see remark \ref{remark:dieumodformpair}) and therefore the isotypical components of $\Hxn$ in the slope decomposition (\ref{eqn:s5slope1}) are isotropic. If we decompose elements $x,y\in \Hxn$ as $x= x_m + x_{m+1}, y= y_m + y_{m+1}$ according to (\ref{eqn:s5slope1}) then $x_m = \io(P(x)), y_m = \io(P(y))$ and $I(x) = I(x_{m+1}), I(y) = I(y_{m+1})$. Together with Lemma \ref{lemma:adj} we obtain
\[
<x,y> = <x_m, y_{m+1}> + <x_{m+1}, y_m> = <P(x), I(y)> - <I(x), P(y)>.
\]
\enddemo

Let $K$ be an extension of $k$ contained in $\C_p$ which is complete with respect to the $p$-adic valuation and of ramification index 1 over $k$ (e.g\ $k=\Q_p$ and $K = \Qur$). We pass from $k$ to $K$ in order to have rational points on $X_{\G}$ (i.e.\ we replace $X_{\G}$ by $(X_{\G})_K$, $V$ by $V_K$ but we will surpress the subscript from now on). The condition on the ramificaton index implies that the image of every $K$-valued point on $\cH_p$ under the reduction map $\red: \cH_p \to |\cT|$ is a vertex). Let $j: U \to X_{\G}$ be the complement of finitely many points $x_1, \ldots, x_r\in X_{\G}(K)$ and choose preimages $z_1, \ldots, z_r\in {\frak F}\seq \cH_p(K)$ of $x_1, \ldots, x_r$ under the projection $\pi: \cH_p \to X_{\G}^{an}$. To simplify the notation we assume that the stabilizer of each point $z_i$ under the action of $\G$ is trivial. 

We have a pairing 
\begin{equation}
\label{eqn:apair7}
< \,\,\,,\,\,\,>_{U} : \Hcu \otimes \Hu \stackrel{\cup}{\lra} H^2_{\DR}(X_{\G}, \cV\otimes \cV) \lra
\end{equation}
 \[
\lra H^2_{\DR,c}(U) \cong  H^2_{\DR}(X_{\G}) \cong k.
\]
The inclusion $j$ induces maps $j_*: \Hcu \to \Hx$, $j^*: \Hx \to \Hu$ such that $< j_*(x),y>_{X_{\G}} = < x, j^*(y)>_{U}$.

We are going to introduce now the analogues of the maps (\ref{eqn:icole}) and (\ref{eqn:isch}) for $\Hcu, \Hu$. We write $\Ind^{\G}(V)$ for the $\G$-module $\Maps(\G, V)$ with $\G$-action is given by $(\ga f)(\tau) \colon = \ga f(\ga^{-1} \tau)$. Let $ad:V \to \Ind^{\G}(V), v \mapsto (\tau \mapsto v)$ and let $\cK^{\bu}(V)$ denote the complex
\[
\begin{CD}
\cK^{\bu}(V) \colon = \Cone( V @> (ad, \ldots, ad)>> \bigoplus_{i=1}^r \Ind^{\G}(V))[-1].
\end{CD}
\]
Let $D$ be the divisor $x_1 + \ldots + x_r$. The de Rham cohomology with compact support $H^i_{\DR,c}(U, \cV)$ can be identified with the hypercohomology
\[
H^i(X_{\G}, [\cL(-D)\otimes \cV \to \Omega^1_{X_{\G}} \otimes \cV]) \cong H^i(\G, \Cone(\Omega^{\bu}(\cH_p)\otimes V \stackrel{\alpha}{\lra} \bigoplus_{i=1}^r \Ind^{\G}(V))[-1])
\]
where $\alpha$ is given (in degree 0) by $F\mapsto (F(z_1), \ldots, F(z_r))$, $F\in \cO(\cH_p)\otimes V$. The distinguished triangle
\[
\cK^{\bu}(V) \lra \Cone(\Omega^{\bu}(\cH_p)\otimes V \stackrel{\alpha}{\lra} \bigoplus_{i=1}^r \Ind^{\G}(V))[-1] \lra C_{\har}(V) \lra \cK^{\bu}(V)[1]
\]
yields short exact sequences
\begin{equation}
\label{eqn:aseq2}
0 \lra H^i(\G, \cK^{\bu}(V)) \stackrel{\io_{U,c}}{\lra} H^i_{\DR,c}(U, \cV)\lra H^{i-1}(\G,C_{\har}(V)) \lra 0.
\end{equation}
and for $i=1$ we get a map $H^1_{\DR,c}(U, \cV)\lra C_{\har}(V)^{\G}$ which we denote by $I_{U,c}$. The injection $\io_{U,c}: H^1(\G, \cK^{\bu}(V)) \to \Hcu$ has again a left inverse $P_{U,c}: \Hcu\to H^1(\G, \cK^{\bu}(V))$ defined in terms of Coleman integration. For that let $\cF(V)$ be the subspace of $V$-valued locally analytic functions on $\cH_p$ which are primitives of elements of $\Omega^1(\cH_p)\otimes V$ (thus the definition of $\cF(V)$ depends on the choice of the branch of the $p$-adic logarithm). Since $V = \Ker(d: \cF(V)\to \Omega^1(\cH_p)\otimes V)$ the complex $[\cF(V) \stackrel{(d, \alpha)}{\lra} \Omega^1(\cH_p)\otimes V\oplus \bigoplus_{i=1}^r \Ind^{\G}(V)]$ (concentrated in degrees 0 and 1) is quasiisomorphic to $\cK^{\bu}(V)$. We defined 
\[
P_{U,c}: \Hcu = H^1(\G, \cO(\cH_p)\otimes V \stackrel{(d, \alpha)}{\lra} \Omega^1(\cH_p)\otimes V\oplus \bigoplus_{i=1}^r \Ind^{\G}(V)]) \lra 
\]
\[
\lra H^1(\G,[\cF(V) \stackrel{(d, \alpha)}{\lra} \Omega^1(\cH_p)\otimes V\oplus \bigoplus_{i=1}^r \Ind^{\G}(V)]) \cong H^1(\G, \cK^{\bu}(V)).
\]
One easily checks that $I_{U,c} = I \circ j_*$ and $q_* \circ P_{U,c} = P \circ j_*$ ($q: V \to \cK^{\bu}(V)$ is the canonical map).

As in section \ref{sec:gltwo} let $C^0(V)$ (resp.\ $C^1(V)$) be the set of maps $f:\cV(\cT) \to V$ (resp.\ maps $f: \edges(\cT) \to V$ such that $f(\bar{e}) = -f(e)$ for all $e\in \edges(\cT)$). We have an exact sequence
\[
0 \lra C_{\har}(V) \lra C^1(V) \stackrel{\delta}{\lra} C^0(V) \lra 0
\]
where $\delta(f)(v) \colon = \sum_{o(e) = v} f(e)$. For $i\in \{1, \ldots, r\}$ we set $v_i = \red(z_i)$ and define $\chi_i: \Ind^{\G}(V) \lra C^0(V)$ by
\[
\chi_i(f)(v)\colon = 
\left\{ \begin{array}{ll}
f(\ga) & \mbox{if $v = \ga v_i$ for some $\ga\in \G$},\\
0 & \mbox{otherwise.}
\end{array}
\right.
\]
Let $C_U(V)$ be the kernel of
\[
C^1(V) \oplus \bigoplus_{i=1}^r \Ind^{\G}(V) \lra C^0(V), (f, f_1, \ldots f_r) \mapsto \delta(f) - \sum_{i=1}^r \, \chi_i(f_i).
\]
We have a short exact sequence
\begin{equation}
\label{eqn:aseq3}
0 \lra C_{\har}(V) \lra C_U(V) \stackrel{\delta}{\lra} \bigoplus_{i=1}^r \Ind^{\G}(V) \lra 0
\end{equation}
For the de Rham cohomology $H^i_{\DR}(U, \cV)$ we have
\[
H^i_{\DR}(U, \cV) = H^i(X_{\G}, \Omega^{\bu}_{X_{\G}}(log(|D|) \otimes \cV]) \cong H^i(\G, \Omega^{\bu}(\cH_p)(log(|\pi^{-1}(D)|)\otimes V).
\]  
There is a distinguished triangle
\[
V \lra \Omega^{\bu}(\cH_p)(log(|\pi^{-1}(D)|)\otimes V \stackrel{I}{\lra} C_U(V)[-1] \lra V[1]
\]
where $I(\omega)(e, \ga_1,  \ldots, \ga_r) \colon = (\Res_e(\omega), \Res_{\ga_1(z_1)}(\omega), \ldots, \Res_{\ga_r(z_r)}(\omega))$. We get short exact sequences
\begin{equation}
\label{eqn:aseq4}
0 \lra H^i(\G, V) \stackrel{\io_{U}}{\lra} H^i_{\DR}(U, \cV)\lra H^{i-1}(\G,C_U(V)) \lra 0.
\end{equation}
For $i=1$ the map $\Hcu\lra C_U(V)^{\G}$ will be denoted by $I_U$.

By identifying $\Hu$ with the space of $\G$-invariant $V$-valued meromorphic differential forms which are of the second kind on $\pi^{-1}(U^{an})$ modulo exact differentials we can define a left inverse $P_U: \Hcu \lra H^1(\G, V)$ of $\io_U$ by the same formula as (\ref{eqn:icole}). We have $P_U \circ j^* = P$ and $I_U \circ j^* = R$.
\begin{lemma}
\label{lemma:adj2}
There exists a pairing
\begin{equation}
\label{eqn:apair8}
< \,\,\,,\,\,\,> = < \,\,\,,\,\,\,>_{\G,U} : H^1(\G, \cK^{\bu}(V))\otimes C_U(V)^{\G} \lra K,
\end{equation}
such that 
\begin{equation}
\label{eqn:adj2}
< \io_{U,c}(x) , y>_U = < x, I_U(y)>_{\G,U}
\end{equation}
for all $x\in H^1(\G, \cK^{\bu}(V))$ and $y\in \Hu$.
\end{lemma}

{\em Proof.} The exact sequences (\ref{eqn:aseq2}), (\ref{eqn:aseq4}) for $i=1$ guarantee the existence of (\ref{eqn:apair8}) with the property (\ref{eqn:adj2}) once  we have shown that $\io_{U,c}(H^1(\G, \cK^{\bu}(V)))$ and  $\io_U(H^1(\G, V))$ are orthogonal under (\ref{eqn:apair7}). This follows by considering the diagram 
\[
\begin{CD}
H^1(\G, \cK^{\bu}(V)) \otimes H^1(\G, V) @>>> H^2(\G, \cK^{\bu}(K))\\
@VV \io \otimes \io V @VV \io V\\
\Hcu \otimes \Hu @>>> H^2_{\DR,c}(U)
\end{CD}
\]
using the fact that $H^2(\G, \cK^{\bu}(K))=0$. The upper map is the composite of the cup-product with the map $H^2(\G,\cK^{\bu}(V) \otimes  V) \to H^2(\G, \cK^{\bu}(K)$ induced by (\ref{eqn:apair1}). 
\enddemo

We also need a concrete description of the pairing (\ref{eqn:apair8}) similar to (\ref{eqn:apair5}). We assume for simplicity that $\G$ is free and leave the formulation of the general case to the reader. Elements of $H^1(\G, \cK^{\bu}(V))$ can be represented by $(r+1)$-tuples $(\fz, f_1, \ldots, f_r)$ such that $ad\circ \fz = \partial(f_1) = \ldots = \partial(f_r)$ where $\fz\in Z^1(\G,V)$, $f_1, \ldots, f_r\in \Ind^{\G}(V)$ and $\partial(f_i)(\ga) = \ga f_i - f_i$. An element of $C_U(V)^{\G}$ is given by a tuple $(g, g_1, \ldots, g_r)\in C^1(V) \oplus\bigoplus_{i=1}^r \Ind^{\G}(V)$ such that $\delta(g) = \sum_{i=1}^r \, \chi_i(g_i)$. With the notation as in (\ref{eqn:apair5}) we have
\begin{equation}
\label{eqn:apair9}
<[(\fz, f_1, \ldots, f_r)], (g, g_1, \ldots, g_r)>_{\G,U}\, = 
\end{equation}
\[
=\sum_{i=1}^g \, <\fz(\ga_i), g(c_i)> + \sum_{j=1,\ldots, r} \, <f_j(1), g_j(1)>.
\]
We leave the verification of this formula as an exercise to the reader.

\begin{conjecture}
\label{conjecture:openshalit}
For all $x\in\Hcu, y\in \Hu$ we have 
\begin{equation}
\label{eqn:ashalit}
<x,y>_U = <P_{U,c}(x),I_U(y)>_{\G, U} - <I_{U,c}(x),P_U(y)>_{\G}.
\end{equation}
\end{conjecture}

It is likely that this can be proved by the methods of (\cite{shalit1}, appendix). Here we will prove it only in the case needed for the application.

\begin{theorem}
Let $X_{\G}$ be the Shimura curve $X$ (over $\Qur$), $V =V_n$ for $n =2m \geq 2$ even and the pairing is pairing {\rm (\ref{eqn:s4pair4})}. Then Conjecture \ref{conjecture:openshalit} holds.
\end{theorem}

{\em Proof.} In this case (\ref{eqn:apair7}) is nondegenerated. Hence there is unique structure of a filtered $(\p,N)$-module on $\Hcu$ such that (\ref{eqn:apair7}) is a map 
\[
\Hcun \otimes \Hun \to \Qur[2m+1]
\]
in $MF_{\Qur}(\p, N)$. By Theorem \ref{theorem:dieumodformgys}, $\Hu$ has a slope decomposition of the form 
\begin{equation}
\label{eqn:a2slope1}
\Hun = \Hun_m \oplus \Hun_{m+1}
\end{equation}
where $\Hun_m = \Hxn_m = \io(H^1(\G,V_n))$ and $\Hxn_{m+1} = \Ker(P_U)$. Thus the slope decompostion of $\Hcun$ is also of the form 
\begin{equation}
\label{eqn:a2slope2}
\Hcun = \Hcun_m \oplus \Hcun_{m+1}
\end{equation}
and $\Hcun_{m+1} = \Hxn_{m+1}$. Since $\io_{U,c}(H^1(\G, \cK^{\bu}(V)))$ and  $\io_U(\break H^1(\G, V))$ are orthogonal we see that $\io_{U,c}(H^1(\G, \cK^{\bu}(V))) \seq \Hcun_m$. But both vector spaces have the same dimension, hence $\io_{U,c}(H^1(\G, \cK^{\bu}(V))) = \break \Hcun_m$. If we write $x\in \Hcun$ as $x= x_m + x_{m+1}$ and $y\in \Hun$ as $y= y_m + y_{m+1}$ according to (\ref{eqn:a2slope1}) and (\ref{eqn:a2slope2}), then $x_m = \io(P_{U,c}(x))$, $y_m = \io_U(P_U(y)) = j^*(\io(P_U(y))$ and we obtain as before
\[
<x,y>_U = <\io(P_{U,c}(x)), y_{m+1}>_U + <j_*(x_{m+1}), \io(P_U(y))>_{X} = 
\]
\[
= <P_{U,c}(x), I_U(y)>_{\G,U} - < I_{U,c}(x), P_U(y)>_{\G}
\]
(the last equality follows from Lemma \ref{lemma:adj} and Lemma \ref{lemma:adj2}). \enddemo

Now assume that $r=2$ and that $\red(z_1) = v = \red(z_2)$. Let $\chi$ be a $\G$-invariant $V_n$-valued meromorphic differential form on $\cH_p$ which is holomorphic outside of $\pi^{-1}(U^{an})$ and which has simple poles at $z_1, z_2$ with $\Res_{z_1}(\chi) = - \Res_{z_2}(\chi)$. The residue theorem implies then that $c_{\chi}(e) \colon = \Res_e(\chi)$ is a harmonic cocycle. We also assume that the class $[\chi]$ lies in the kernel of $N: \Hun\to\Hun$. Note that this implies that $c_{\chi} = 0$ (see \cite{shalit1}). 

\begin{coro} 
\label{coro:openshalit}
For every $f\in M_k(\G)$ we have
\[
< I([\omega_f)], P_U(\chi)>_{\G} = < F_{\omega_f}(z_1), \Res_{z_1}(\chi)> + 
< F_{\omega_f}(z_2), \Res_{z_2}(\chi)>.
\]
\end{coro}

{\em Proof.} Since $F^{m+1} \Hcun$ and $F^{m+1} \Hun$ are orthogonal we have $<[\omega_f], [\chi]>_U =0$. Hence by (\ref{eqn:ashalit}) we get $< I([\omega_f)], P_U(\chi)>_{\G} = < P_{U,c}([\omega_f]), I_U([\chi]>_{\G,U}$. However using the formula (\ref{eqn:apair9}) it is easy to see that
\[
< P_{U,c}([\omega_f]), I_U([\chi])>_{\G,U} =
\]
\[
= < F_{\omega_f}(z_1), \Res_{z_1}(\chi)> + < F_{\omega_f}(z_2), \Res_{z_2}(\chi)>.
\]
\enddemo

\bigskip
\begin{tabbing}
 \hspace{1cm}\= Adrian Iovita\hspace{3cm}       \= Michael Spiess \\
 \>Department of Mathematics \> School of Mathematical Sciences \\
 \>University of Washington \> University of Nottingham \\
  \>Box 355754 \> University Park \\ 
 \>Seattle, WA 98105 \> NG7 2RD Nottingham \\ 
 \>USA \> United Kingdom 
 \end{tabbing}

\noindent \hspace{.9cm} e-mail: iovita@math.washington.edu \\ 
\hspace*{2.3cm} mks@maths.nott.ac.uk 


\begin{thebibliography}{XXXX}

\bibitem[Ber]{bert} P.\ Berthelot, {\em Cohomologie rigide et cohomologie rigid \`a supports propre},

\bibitem[Be]{besser} A.\ Besser, {\em CM cycles over Shimura curves}, Journal of Algebraic Geometry {\bf 4} (1994).

\bibitem[BD1]{bertolini_darmon1}
M.\ Bertolini and H.\ Darmon, {\em Heegner points on Mumford-Tate curves}, Invent.\ Math.\ {\bf 126} (1996), 413--456.

\bibitem[BD2]{bertolini_darmon2}
M.\ Bertolini and H.\ Darmon, {\em Heegner points, $p$-adic $L$-functions
and the Cerednik-Drinfeld uniformisation}, Invent.\ Math.\ {\bf 131} (1998), 453--491.

\bibitem[BDIS]{bdis}
M.\ Bertolini, H.\ Darmon, A.\ Iovita and M.\ Spie{\ss}, {\em Teitelbaum's exceptional zero conjecture in the anticyclotomic setting}, to appear in Amer.\ J.\ Math.\

\bibitem[BK]{bloch_kato}
S.\ Bloch and K.\ Kato, {\em $L$-functions and Tamagawa numbers of motives.} In: The Grothendieck Festschrift I, 333--400, Progr.\ Math.\ {\bf 108}, Birkh{\"a}user, Boston (1993).

\bibitem[BC]{boutot_carayol}
J.\-F.\ Boutot and H.\ Carayol, {\em Uniformisation $p$-adique des courbes de Shimura: les th\'eor\`emes de \v Cerednik et de Drinfeld.}  Courbes modulaires et courbes de Shimura (Orsay, 1987/1988). Ast\'erisque {\bf 196-197} (1991), {\bf 7}, 45--158 (1992).

\bibitem[Ce]{cerednik}
I.V.\ ~\v Cerednik, {\em Uniformisation of algebraic curves by discrete
arithmetic subgroups of ${\rm PGL}\sb{2}(k\sb{w})$ with compact quotient
spaces.} (Russian) Mat.\ Sb.\ (N.S.) {\bf 100} (1976), 59--88.

\bibitem[Co]{cole} R.\ Coleman, {\em A $p$-adic Shimura Isomorphism and $p$-adic Periods of Modular forms.} In: $p$-Adic Monodromy and the Birch and Swinnerton-Dyer Conjecture (Boston, 1991). Contemporary Math.\ {\bf 165} (1994), 21--51.

\bibitem[CI]{coleman_iovita} R.\ Coleman and A.\ Iovita, {\em Frobenius and Monodromy operators on semistable curves via degeneration}, preprint.

\bibitem[De1]{deligne1} P.\ Deligne, {\em Travaux de Shimura}, 
S\'eminaire Bourbaki (1970/71), Exp.\ No.\ 389. Lect.\ Notes in Math.\ {\bf 244}, Springer Verlag 1971, 123--165.

\bibitem[De2]{deligne2}
P.\ Deligne, {\em Hodge cycles on abelian varieties}. In: Hodge Cycles, Motives and Shimura Varieties (P.\ Deligne, J.\ Milne, A.\ Ogus, K.\ Shih). Lect.\ Notes in Math.\ {\bf 900}, Springer Verlag 1982, 9--100.

\bibitem[DM]{de-mu} C.\ Deninger and J.\ Murre, {\em Motivic decomposition of abelian schemes and the Fourier transform}, J.\ reine angew.\ Math.\ {\bf 422} (1991), 201--219.

\bibitem[Dr]{drinfeld} V.G.\ Drinfeld, {\em Coverings of $p$-adic symmetric domains.} (Russian) Funkcional. Anal.\ i Prilo\v zen.\ {\bf 10} (1976), 29--40.

\bibitem[Fa1]{faltings1} G. Faltings, {\em Crystalline cohomology and $p$-adic Galois-representations}. In: Algebraic analysis, geometry, and number theory (Baltimore, MD, 1988). Johns Hopkins Univ.\ Press 1989, 25--80.

\bibitem[Fa2]{faltings2} G. Faltings, {\em Cristalline cohomology of semistable curve ---the $\Q_p$-theory}, J.\ Algebraic Geometry {\bf 6} (1997), 1--18.

\bibitem[Fo1]{font1} J.\-M.\ Fontaine, {\em Les corps des p\'eriodes $p$-adiques}, P\'eriodes $p$-adiques. S\'eminaire de Bures (1988). Ast\'erisque {\bf 223} (1994), 59--101.

\bibitem[Fo2]{font2} J.\-M.\ Fontaine, {\em Repr\'esentations $p$-adiques semi-stables}, P\'eriodes $p$-adiques. S\'eminaire de Bures (1988). Ast\'erisque {\bf 223} (1994), 113--184.

\bibitem[Hy]{hyodo} O.\ Hyodo, {\em $H^1_g(K, V) = H^1_{\st}(K, V)$}, unpublished manuscript.

\bibitem[Ja]{jannsen} U.\ Jannsen, {\em Continuous {\'E}tale Cohomology}, Math.\ Ann.\ {\bf 280} (1988), 207--245.

\bibitem[Ka]{kato} K.\ Kato, {\em Logarithmic structures of Fontaine-Illusie}. In: Algebraic analysis, geometry, and number theory (Baltimore, MD, 1988). Johns Hopkins Univ.\ Press 1989, 191--224.

\bibitem[Ki]{kisin} M.\ Kisin, {\em Potential semi-stability of $p$-adic \'etale cohomology}, to appear in Israel J.\ of Math.\

\bibitem[K{\"u}]{kuennemann} K.\ K{\"u}nnemann, {\em On the Chow Motive of an Abelian Scheme.} In: Motives (Seattle, WA, 1991). Proceedings of Symposia in Pure Math. {\bf 55.1} (1994), 189--205. 

\bibitem[Ma]{mazur} B.\ Mazur, {\em On monodromy invariants occurring in global arithmetic and Fontaine's theory.} In: $p$-Adic Monodromy and the Birch and Swinnerton-Dyer Conjecture (Boston, 1991). Contemporary Math. {\bf 165} (1994), 1--20.

\bibitem[MTT]{mazur_tate_teitelbaum}
B.\ Mazur, J.\ Tate, and J.\ Teitelbaum,
{\em On $p$-adic analogues of the conjectures of Birch and
Swinnerton-Dyer}, Invent.\ Math.\ {\bf 84} (1986), 1--48.

\bibitem[Ne1]{nekovar1} J.\ Nekovar, {\em Kolyvagin's method for Chow groups of Kuga-Sato varieties}. Invent.\ Math.\ {\bf 107} (1992), 99--125.

\bibitem[Ne2]{nekovar2} J.\ Nekovar, {\em On $p$-adic height pairings}, S\'eminaire de Th\'eorie des Nombres, Paris, 1990--91, 127--202, Progr.\ Math.\ {\bf 108}, Birkh{\"a}user, Boston (1993).
 
\bibitem[Ne3]{nekovar3} J.\ Nekovar, {\em On the $p$-adic height of Heegner cycles}, Math.\ Ann.\ {\bf 302} (1995), 609--686. 

\bibitem[Ne4]{nekovar4} J.\ Nekovar, {\em $p$-adic Abel-Jacobi maps and $p$-adic heights.} In: The arithmetic and geometry of algebraic cycles (Banff, AB, 1998), CRM Proc.\ Lecture Notes, 24, Amer.\ Math.\ Soc.\, Providence, RI (2000), 367--379.

\bibitem[Ne5]{nekovar5} J.\ Nekovar, {\em Syntomic cohomology and $p$-adic regulators}, preprint (1998).

\bibitem[Og]{ogus} A.\ Ogus, {\em $F$-Isocrystals and de Rham Cohomology II --
Convergent Isocrystals}, Duke Math.\ J.\ {\bf 51} (1984), 765--850.

\bibitem[Ri1]{ribet1} K.\ Ribet, {\em Sur les vari\'et\'es ab\'eliennes \`a multiplications r\'eelles}, C.\ R.\ Acad.\ Sc.\ Paris {\bf 291} (1980), 121--123.

\bibitem[Ri2]{ribet2} K.\ Ribet, {\em On $l$-adic representations attached to modular forms. II}, Glasgow Math.\ J.\ {\bf 27} (1985), 185--194.

\bibitem[RZ]{rap-zink} M.\ Rapoport and T.\ Zink, {\em Period spaces for $p$-divisible groups}. Annals of Mathematics Studies {\bf 141}, Princeton University Press (1996).

\bibitem[Sch1]{scholl} A.J.\ Scholl, {\em Motives for modular forms}, Invent.\ Math.\ {\bf 100} (1990), 419--430.

\bibitem[Sch2]{scholl2} A.J.\ Scholl, {\em Classical motives.} In: Motives (Seattle, WA, 1991). Proceedings of Symposia in Pure Math. {\bf 55.1} (1994), 163--187. 

\bibitem[dS1]{shalit1} E.\ de Shalit, {\em Eichler cohomology and periods of modular forms on $p$-adic Schottky groups}, J.\ reine angew. Math. {\bf 400} (1989), 3--31.

\bibitem[dS2]{shalit2} E.\ de Shalit, {\em A formula for the cup product on Mumford curves}, S{\'e}minaire de Th{\'e}orie des Nombres, 1987--1988 (Talence, 1987--1988), Exp.\ No.\ 47, 10 pp.\, Univ.\ Bordeaux I.

\bibitem[dS3]{shalit3} E.\ de Shalit, {\em Differentials of the second kind on Mumford curves}, Israel Journal of Mathematics {\bf 71} (1990), 1--16.

\bibitem[Ts]{tsuji} T.\ Tsuji, {\em $p$-adic {\'e}tale cohomology and crystalline in the semi-stable reduction case}, Invent.\ Math.\ {\bf 137} (1999), 233--411.

\bibitem[Te]{teitelbaum} J. Teitelbaum, {\em Values of $p$-adic $L$-functions and a $p$-adic Poisson kernel}, Invent. Math. {\bf 101} (1990), 395--410.

\bibitem[Wo]{wortmann} S.\ Wortmann, work in progress.

\bibitem[Zi]{zink} T.\ Zink, {\em Catiertheorie commutativer formaler Gruppen}, Teubner Texte zur Mathematik {\bf 68}, Leipzig 1984.


\end{thebibliography}
\end{document}